\input amstex
\documentstyle{amsppt}
\magnification=\magstep1

\NoBlackBoxes
\TagsAsMath

\pagewidth{6.5truein}
\pageheight{9.0truein}

\loadbold

\long\def\ignore#1\endignore{#1}

\ignore
\input xy \xyoption{matrix} \xyoption{arrow}
          \xyoption{curve}  \xyoption{frame}
\def\edge{\ar@{-}}
\def\dttdar{\ar@{.>}}

\def\levelpool#1{\save [0,0]+(-3,3);[0,#1]+(3,-3)
                  **\frm<10pt>{.}\restore}
\def\dashedge{\ar@{--}}

\def\dropdown#1{\save+<0ex,-4ex> \drop{#1} \restore}
\def\dropup#1{\save+<0ex,4ex> \drop{#1} \restore}
\endignore

\def\seq{\mathrel{\widehat{=}}}
\def\la{{\Lambda}}
\def\lamod{\Lambda\text{-}\roman{mod}}
\def\Lamod{\Lambda\text{-}\roman{Mod}}
\def \len{\operatorname{length}} 
\def\AA{{\Bbb A}}
\def\CC{{\Bbb C}}
\def\PP{{\Bbb P}}
\def\SS{{\Bbb S}}

\def\NN{{\Bbb N}}

\def\aut{\operatorname{Aut}}
\def\Aut{\operatorname{Aut}}

\def\stab{\operatorname{Stab}}

\def\GL{\operatorname{GL}}

\def\Ker{\operatorname{Ker}}
\def\Im{\operatorname{Im}}

\def\autlap{\operatorname{Aut}_\la(P)}
\def\bigautlap{\operatorname{Aut}_\la({\bold P})}

\def\End{\operatorname{End}}

\def\can{\operatorname{can}}

\def\A{{\Cal A}}

\def\fraka{{\frak a}} 
\def\frakb{{\frak b}} 
 
\def\C{{\frak C}}

\def\bd{{\bold{d}}}

\def\I{{\Cal I}}

\def\T{{\Cal T}}

\def\m{{\frak m}}

\def\P{{\widehat{\frak P}}}
\def\bp{\bold P}
\def\bphat{\bold{\widehat{P}}}
\def\bphi{\bold{{\Phi}}}
\def\bpsi{\bold{{\Psi}}}

\def\S{{\sigma}}

\def\T{{\Cal T}}

\def\bt{{\bold{t}}}

\def\etahat{\widehat{\eta}}

\def\Jhat{\Cal J}
\def\Phat{\widehat{P}}
\def\Qhat{\widehat{Q}}
\def\Rhat{\widehat{R}}

\def\zhat{\widehat{z}}

\def\GRASS{\operatorname{\text{\smc{grass}}}}

\def\modlad{\operatorname{\bold{Mod}}_d(\Lambda)}
\def\modlabd{\operatorname{\bold{Mod_d}}(\Lambda)}

\def\toptd{\operatorname{\bold{Mod}}^T_d}
\def\toptbd{\operatorname{\bold{Mod}}^T_{\bold d}}

\def\laySS{\operatorname{\bold{Mod}}(\SS)}
\def\lay#1{\operatorname{\bold{Mod}}(#1)}

\def\Hom{\operatorname{Hom}}

\def\Schu{\operatorname{\Schu}}
\def\bigSchu{\operatorname{\ssize{SCHU}}}

\def\bigSchuS{\operatorname{\ssize{SCHU}}(\sigma)}
\def\grasstd{\operatorname{\frak{Grass}}^T_d}
\def\grasstbd{\operatorname{\frak{Grass}}^T_{\bold d}}
\def\biggrasstd{\GRASS^T_d}
\def\biggrasstbd{\GRASS^T_{\bold d}}
\def\grassd{\GRASS_d(\Lambda)}
\def\grassbd{\GRASS_{\bold d}(\Lambda)}
\def\Grfrak{\operatorname{\frak{Gr}}}

\def\autlap{\aut_\Lambda(P)}

\def\unirad{\bigl(\aut_\la(P)\bigr)_u}
\def\bigunirad{\bigl(\aut_\la(\bp)\bigr)_u}
\def\grassS{\operatorname{\frak{Grass}}(\S)}

\def\grassSS{\operatorname{\frak{Grass}}(\SS)}
\def\biggrassSS{\GRASS(\SS)}

\def\biggrassS{\GRASS(\S)}

\def\grass{\operatorname{\frak{Grass}}}
\def\biggrass{\GRASS}

\def\degen{\le_{\text{deg}}}
\def\underbardim{\operatorname{\underline{dim}}}
\def\id{\operatorname{id}}

\def\Schu{\operatorname{Schu}}

\def\vsubseteq{\hbox{$\bigcup$\kern0.1em\raise0.05ex\hbox{$\tsize|$}}}

\def\generic{{\bf 1}}
\def\codes{{\bf 2}}
\def\GeomII{{\bf 3}}
\def\GeomIV{{\bf 4}}
\def\Bor{{\bf 5}}
\def\Gab{{\bf 6}}
\def\menace{{\bf 7}}
\def\grassI{{\bf 8}} 
\def\degen{{\bf 9}}
\def\degenII{{\bf 10}}
\def\Ros{{\bf 11}}
\def\RosII{{\bf 12}}
\def\Zwara{{\bf 13}}

\topmatter

\title A hierarchy of parametrizing varieties
for representations
\endtitle


\author  B. Huisgen-Zimmermann
\endauthor

\address Department of Mathematics, University of California, Santa
Barbara, CA 93106 \endaddress

\email birge\@math.ucsb.edu \endemail

\dedicatory Dedicated to Carl Faith and Barbara Osofsky \enddedicatory

\thanks The author was partly supported by a grant
from the National Science Foundation. 
\endthanks

\abstract The primary purpose is to introduce and explore projective
varieties, $\grassbd$, parametrizing the  full collection of those
modules over a finite dimensional algebra $\la$ which have dimension
vector $\bd$. These varieties extend the smaller varieties previously
studied by the author; namely, the projective varieties encoding those
modules with dimension vector $\bd$ which, in  addition, have a preassigned
top or radical layering. Each of the $\grassbd$ is again partitioned by
the action of a linear algebraic group, and covered by certain
representation-theoretically defined affine subvarieties which are stable
under the unipotent radical of the acting group. A special case of the
pertinent theorem served as a cornerstone in the
work on generic representations by Babson, Thomas, and the author.
Moreover, applications are given to the study of degenerations.
\endabstract

\subjclassyear{2000}
\subjclass 16G10, 16G20, 14L30 \endsubjclass


\endtopmatter

\document

\head 1. Introduction and notation \endhead

Our primary aim is to extend some of the concepts, constructions and
results from
\cite{\grassI}, \cite{\degen}, and \cite{\generic} for wider
applicability towards exploring the representation theory of a basic
finite dimensional algebra $\la$ over an algebraically closed field $K$.
These varieties constitute the foundation for present and future work on
degenerations, irreducible components and generic representations. 
More specifically, our purpose is to introduce a projective parametrizing
variety $\grassbd$ for the full collection of isomorphism classes of
modules with a given dimension vector $\bd = (d_1, \dots, d_n)$ recording
the composition multiplicities of the simple objects $S_1, \dots, S_n$ in
$\lamod$.  This large variety is to supplement the projective variety
$\grasstbd$ and its subvariety $\grassSS$, which encode the
$\bd$-dimensional modules $M$ with fixed top
$T = M/JM$ and fixed radical layering $\SS = (J^l M/ J^{l+1}M)_{0
\le l \le L}$, respectively \cite{\grassI, \degen, \generic}; here
$J$ denotes the Jacobson radical of $\la$ and $L+1$ its nilpotency
index.  This
permits us to study ``global" problems by moving back and forth between
local and global settings.  (``$\biggrass$" and ``$\grass$" both stand for
``Grassmann", indicating that the varieties carrying one of these labels
are subvarieties of Grassmannians.)    

The new variety
$\grassbd$ and its subvarieties $\biggrasstbd$ and $\biggrassSS$  -- the
latter parametrizing the same isomorphism classes of modules as $\grasstd$
and $\grassSS$, respectively  --  again come equipped with an algebraic
group action; as in all of the previously considered cases, the orbits
partition the given variety into strata representing the different
isomorphism classes of modules.  The acting group in the new setting is in
turn larger than the one that accompanies the smaller varieties $\grasstd$
and $\grassSS$, whence the enlarged projective setting $\biggrass$ loses
the advantage of comparatively small orbit dimensions.  However, it
retains several crucial assets, notably completeness.   An additional
plus of all the varieties arising in the Grassmannian scenario lies in the
fact that the acting group typically has a large unipotent radical; this
brings a cache of methods to bear which supplement those available for the
reductive group actions in the affine scenario.  Yet, more importantly,
the Grassmannian varieties $\biggrassSS$ and
$\grassSS$ are covered by certain representation-theoretically  defined
open subvarieties -- $\biggrassS$ and $\grassS$ respectively  --  which are
particularly accessible:  Namely,
the modules with fixed radical layering $\SS$ are further subdivided into
classes of modules characterized by much finer structural invariants,
their ``skeleta" $\sigma$ (see Definition 3.1).  On the level of the
$\grassS$, representation-theoretic investigations can be carried out more
effectively than on the levels located above, with good venues for moving
information back ``up the ladder".   

To benefit from the increased representation-theoretic transparency of the
 variety
$\grasstbd$ representing the same objects as $\biggrasstbd$, the 
practical philosophy is as follows:  When a global problem tackled in
the large variety $\grassbd$ has been played down to $\biggrasstbd$ 
(non-closed and hence non-projective), consisting of the points encoding
modules with a fixed top $T$, it is advantageous to bring the smaller
projective counterpart $\grasstbd$ into the picture.  The governing idea
is to reduce any given problem to the smallest possible setting to render
it more tractable and then move in reverse.  So, in other
words, it is advantageous to explore both the large and small scenarios
side by side.      
     
In order to place
the large varieties $\grassbd$ into context, we begin by reminding the
reader of the other candidates on the previously defined list of
nested ``small" parametrizing varieties of representations:   
$$\grasstd \supseteq \grasstbd \supseteq \grassSS \supseteq \grassS.$$
Here $d$ is the total dimension of $\bd$, and $\grasstd$ parametrizes 
the modules with dimension $d$ and top $T$.  By $\SS = (\SS_0, \dots,
\SS_L)$ we denote a sequence of semisimple modules with ``top" $\SS_0 =
T$ such that $\bd = \underline{\dim} \bigoplus_{0 \le l
\le L} \SS_l$; it represents a typical radical layering of a module with
top $T$ and dimension vector $\bd$.  In order to descend an additional step
to the varieties labeled 
$\grassS$, we present $\la$ as a path
algebra modulo relations, that is, in the form
$\la = K\Gamma/I$, where $\Gamma$ is the Gabriel quiver of $\la$, and $I$
an admissible ideal in the path algebra $K\Gamma$.  This allows us to
organize
$\la$-modules in terms of their ``skeleta", roughly speaking, path bases
of $\la$-modules which reflect their underlying $K\Gamma$-module
structure.  Starting with a skeleton $\S$ compatible with the sequence
$\SS$ (Definition 3.5), we denote by
$\grassS$ the subvariety of $\grassSS$ which consists of the points
corresponding to modules with skeleton $\S$.  

We will keep
the old notation for compatibility with prior work, and distinguish the
larger varieties to be defined and discussed in the sequel by capital
letters in order to place emphasis on the big versus the small scenario. 
The enlarged setting consists of 
$$\grassd \supseteq \grassbd \supseteq \biggrasstbd \supseteq \biggrassSS
\supset \biggrassS,$$ 
where $\grassd$ and $\grassbd$ are projective varieties representing all
modules with dimension $d$, resp\. dimension vector
$\bd$, and $\biggrasstbd$, $\biggrassSS$, $\biggrassS$ are the locally
closed subvarieties consisting of the points which, in analogy with the
previous list, correspond to the modules with top
$T$, radical layering $\SS$, and skeleton $\S$, respectively.  
\bigskip

We preview the main results of the paper, Theorems 3.12 and 3.17, in a
loosely phrased version.  The ``small version" has
already been used without proof, with reference to the present article, in
\cite{\generic} and \cite{\degen}; in both of these investigations, it
plays a pivotal role.  Moreover, it served as basis for an algorithm,
developed by  Babson, Thomas and the author \cite{\codes}, which computes
polynomials for the smallest of the considered varieties  --  that is, for
the 
$\grassS$  --  from the quiver $\Gamma$ and generators for the ideal $I$.  

Suppose $P$ and $\bp$ are  projective covers of $T$ and $\bigoplus_{1 \le
i \le n} S_i^{d_i}$, respectively.  Then $\grasstbd$ and $\biggrasstbd$
consist of the submodules $C$ of $JP$, resp. $\bp$, with the property
that the factor module
$P/C$, resp. $\bp / C$, has top $T$ and dimension vector $\bd$.  Denote
by $\phi$, resp. $\bphi$, the map sending
$C$ to the isomorphism class of $P/C$, resp.
$\bp/C$.

\proclaim{Theorem} Let
$\S$ be a skeleton of dimension vector $\bd$ which is compatible
with $\SS$.  Then the set $\grassS$ consisting of
the points in $\grasstbd$ which correspond to modules with
skeleton $\S$ is an open subvariety of $\grassSS$; analogously,
$\biggrassS$ is an open subvariety of $\biggrassSS$.

There is an isomorphism $\psi$ from $\grassS$ onto a closed subvariety of
an affine space
$\AA^N$, together with an explicitly specified map $\chi$ from $\Im(\psi)$
onto the set of isomorphism classes of modules with skeleton $\S$, such
that the following triangle commutes:

$$\xymatrixrowsep{2pc}\xymatrixcolsep{4pc}
\xymatrix{
\Im(\psi) \save+<-6.6ex,0.2ex> \drop{\AA^N\supseteq} \restore
\ar[dd]_{\psi^{-1}} \ar[dr]^-{\chi}\\
 &\save+<23ex,0ex> \drop{\txt{ \{\rm isom.~types of
$\la$-modules with skeleton $\S$\} }} \restore\\
\grassS \save+<-10.1ex,0.2ex> \drop{\grasstd\supseteq} \restore
\ar[ur]_-{\phi} }.$$

Polynomial equations for the isomorphic copy $\Im(\psi)$ of
$\grassS$ in $\AA^N$ can be algorithmically obtained from
$\Gamma$ and left ideal generators for $I$. 

In the big scenario, there exist an isomorphism
$\bpsi$ from $\biggrassS$ to a closed subset of an affine
space and a map $\boldsymbol{\chi}$ from $\Im(\bpsi)$ to the set
of isomorphism classes of $\la$-modules with skeleton $\S$, giving rise to
a corresponding commutative diagram.

Moreover, $\biggrassS \cong \grassS \times \AA^{N_0}$, where $N_0$ is
determined by the semisimple sequence $\SS$ alone.
\endproclaim

We follow with a schematic overview of the relevant varieties.  Those in
the first three slots of the left column below carry a canonical action by
the big automorphism group
$\bigautlap$, those in the corresponding slots of the central column are
equipped with the conjugation action by
$\GL_d$, and the first two varieties in the rightmost column with an action
by $\autlap$ that parallels that of $\bigautlap$.  Moreover, the horizontal
double arrows indicate easy transfer of information.  
The varieties on the same level parametrize the
same set of isomorphism types of modules, with geometric information being
transferrable between the two sides via Propositions 2.1 and 2.5. 
In the
first three rows, this concerns information about action-stable subsets. 
The affine varieties in the last row are not stable under $\bigautlap$ or
$\autlap$ in general, but only stable under the
action of the unipotent radical of $\bigautlap$, resp. of
$\autlap$, and a maximal torus in the relevant automorphism group; on this
lowest level, the transfer of information is by way of the final statement
of the above theorem.

$$\xymatrixrowsep{0.5pc}\xymatrixcolsep{4pc}
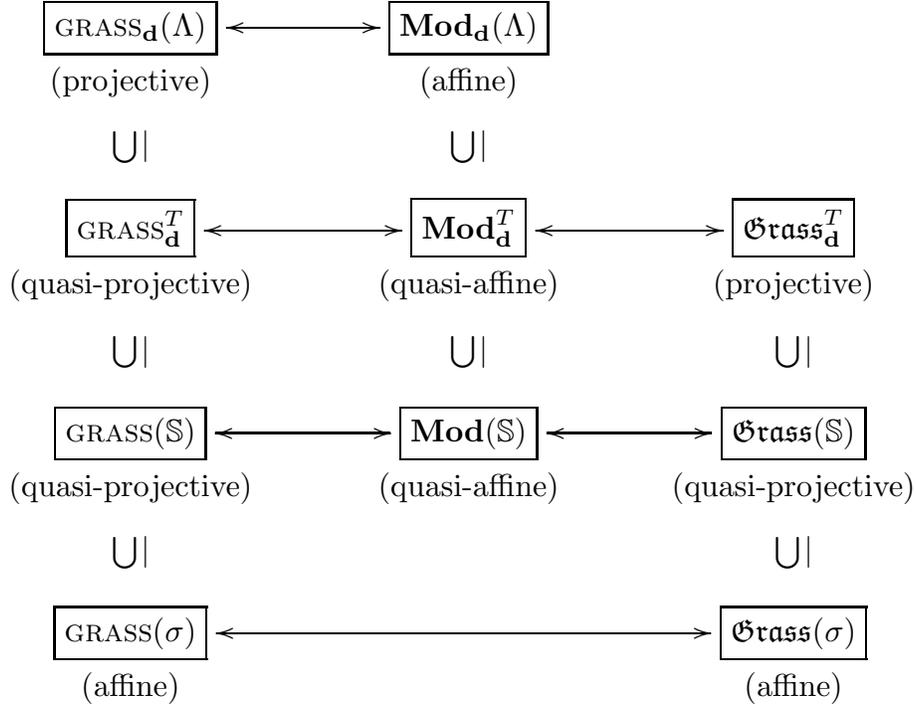
\xymatrix{
\boxed{\grassbd} \dropdown{\txt{(projective)}} \ar@{<->}[r] &\boxed{\modlabd}
\dropdown{\txt{(affine)}} \\  \\
\vsubseteq &\vsubseteq \\
\boxed{\biggrasstbd} \dropdown{\txt{(quasi-projective)}} \ar@{<->}[r]
&\boxed{\toptbd} \dropdown{\txt{(quasi-affine)}}
\ar@{<->}[r] &\boxed{\grasstbd} \dropdown{\txt{(projective)}} \\  \\
\vsubseteq &\vsubseteq &\vsubseteq \\
\boxed{\biggrassSS} \dropdown{\txt{(quasi-projective)}} \ar@{<->}[r]
&\boxed{\laySS} \dropdown{\txt{(quasi-affine)}}
\ar@{<->}[r] &\boxed{\grassSS} \dropdown{\txt{(quasi-projective)}} \\  \\
\vsubseteq &&\vsubseteq \\
\boxed{\biggrassS} \dropdown{\txt{(affine)}} \ar@{<->}[rr] &&\boxed{\grassS}
\dropdown{\txt{(affine)}}
 }$$
\bigskip    

\noindent In the final section, we present some first applications of the
new Grassmannians to degenerations.
\medskip

Throughout, $\la$ denotes a basic finite dimensional algebra
over an algebraically closed field $K$.  Hence,
we do not lose generality in assuming that $\la = K\Gamma/I$ is a path
algebra modulo relations as above.   
The vertices $e_1, \dots, e_n$ of $\Gamma$ will be
identified with the primitive idempotents of $\la$ corresponding
to  the paths of length zero.  As is well-known, the left ideals $\la e_i$
then represent all indecomposable projective (left) $\la$-modules, up
to isomorphism, and the factors $S_i = \la e_i /J e_i$, where $J$ is
the Jacobson radical of $\la$, form a set of representatives for the simple
left $\la$-modules.  By
$L+1$ we will denote the Loewy length of $\la$.   Moreover, we will observe
the following conventions:  The product $pq$ of two paths $p$ and $q$ in
$K\Gamma$ stands for ``first $q$, then $p$"; in particular, $pq$ is zero
unless the end point of $q$ coincides with the starting point of $p$.  In
accordance with this convention, we call a path $p_1$ an {\it initial
subpath\/} of $p$ if
$p = p_2 p_1$ for some path $p_2$.  A {\it path in $\la$\/} is a residue
class of the form $p + I$, where
$p$ is a path in $K Q \setminus I$; we will suppress the residue
notation, provided there is no risk of ambiguity.  Further, we will
gloss over the distinction between the left $\la$-structure of a module
$M \in \Lamod$ and its induced $K\Gamma$-module structure when there is no
danger of confusion.  An element
$x$ of $M$ will be called a {\it top element\/} of
$M$ if $x \notin JM$ and $x$ is {\it normed\/} by some $e_i$, meaning
that $x = e_i x$.  Any collection $x_1, \dots, x_m$ of top elements of $M$
generating
$M$ and linearly independent modulo $JM$ will be referred to as a {\it full
sequence of top elements in $M$\/}.

\head 2.  Old and new parametrizing varieties for representations
\endhead  

We begin by laying out the various varieties listed in the three upper
rows of the above diagram. In subsection 2.A, we briefly remind the reader
of the classical affine parametrizing varieties
$\modlabd$. In 2.B, we review crucial properties of the subvariety
$\toptbd$ and its projective counterpart $\grasstbd$.  In 2.C, we cut down
further in size to $\laySS$ and $\grassSS$.  Next, in 2.D, we introduce the
larger projective variety $\grassbd$ and its subvarieties $\biggrasstbd$
and $\biggrassSS$, and compare the ``big'' and ``small'' settings
following the horizontal double arrows.  In 2.E, we use the sets $\laySS$
and
$\biggrassSS$ to build useful closed subvarieties of $\modlabd$ and
$\grassbd$. In 2.F, finally, we
discuss the structure of the automorphism groups acting on the projective
parametrizing varieties.

\subhead 2.A. Reminder: The classical affine variety 
$\modlabd$
\endsubhead
 
Let $a_1, \dots, a_r$ be a set of algebra generators for
$\la$ over $K$.  A convenient set
of such generators consists of the primitive idempotents (= vertices)
$e_1, \dots, e_n$ together with the (residue classes in $\la$ of the)
arrows in $\Gamma$.  Recall that, for any  natural number $d$, the classical
affine variety of $d$-dimensional representations of $\la$ can be
described in the form
$$\multline
\modlad =  \\
 \{(x_i) \in \prod_{1 \le i \le r} \End_K(K^d)
\mid \text{\ the\ } x_i \text{\ satisfy all relations satisfied by the
\ } a_i\}.
\endmultline$$
  As is well-known, the isomorphism classes of $d$-dimensional (left)
$\la$-modules are in one-to-one correspondence with the orbits of
$\modlad$ under the $\GL_d$-conjugation action.  Moreover, the
connected components of $\modlad$ are known to be in natural one-to-one
correspondence with the dimension vectors $\bd = (d_1, \dots, d_n)$ of
total dimension $d$, meaning $|\bd|= \sum_i d_i = d$.  Namely, the
connected component corresponding to $\bd$ is
$$\multline
\modlabd =  \\
\{x \in \modlad \mid \text{the module corresponding to}\ x
\ \text{has dimension vector}\ \bd \};
\endmultline$$
see \cite{\Gab}, Corollary 1.4.
In particular, $\modlabd$ is a closed subvariety of $\modlad$, and
hence again affine.  In the following, we will exclusively focus on
varieties parametrizing modules with fixed dimension vector $\bd$ (as
opposed to fixed dimension $d$), since that will not result in any
restrictions to applicability.    

If $I = 0$, that is, if $\la$ is hereditary, the
connected varieties $\modlabd$ are even irreducible, but
this fails already in small non-hereditary examples.  For instance, when
$\bd = (1,1)$ and  $\la = K\Gamma/I$, where $\Gamma$ is the quiver  
\ignore $\xymatrixrowsep{1pc}\xymatrixcolsep{2pc} \xymatrix{  
1 \ar@/^0.2pc/[r] &2 \ar@/^0.2pc/[l] }$\endignore
and $I$ is generated by the paths of length $2$, the variety $\modlabd$,
and analogously $\grassbd$, has two irreducible components, namely the
orbit closures of the two uniserial modules of dimension $2$ (see Examples
2.8(1) for more detail in the Grassmannian setting). 

\subhead 2.B. The quasi-affine variety
$\toptbd$ and its counterpart, the projective variety
$\grasstbd$
\endsubhead  

The two isomorphism invariants of a $\la$-module $M$ which will be
pivotal here are the {\it top\/} and the {\it radical
layering\/} of $M$, the latter being a refinement of the former. 
The {\it top\/} of $M$ is defined as
$M/JM$,  and the {\it radical layering\/} as the sequence
$\SS(M) = \bigl(J^lM / J^{l+1}M \bigr)_{0\le l \le L}$ of semisimple
modules.  We will say that a  module $M$ has {\it top $T$\/} in case $M/JM
\cong T$. In fact, we will identify  isomorphic semisimple
modules, and thus write $M/JM  = T$ in this situation.  In
light of our identification, the set of all finite dimensional semisimple
modules is partially ordered under inclusion. Specifically, we write $T
\le T'$ to denote that a semisimple module $T$ is (isomorphic to) a
submodule of a semisimple module $T'$.
\medskip

\noindent {\it The following choices and notation will be observed
throughout\/}: We fix a semisimple $\la$-mod\-ule $T$, say 
$$T = \bigoplus_{1 \le i \le n} S_i^{t_i}\,,$$ 
with dimension vector $\bt = (t_1, \dots, t_n)$ and total dimension $t
=\sum_i t_i$, and denote by
$$P = \bigoplus_{1 \le r \le t} \la z_r$$
the {\it distinguished projective cover\/} of $T$; this means that $\la
z_1,\dots,\la z_{t_1}$ are isomorphic to $\la e_1$, while $\la
z_{t_1+1},\dots,\la z_{t_1+t_2}$ are isomorphic to $\la e_2$, and so on.
Here each $z_r$ is a top element of $P$ normed by a primitive idempotent
$e(r)$, that is,
$z_r = e(r) z_r$, and we refer to $z_1, \dots, z_t$
as the {\it distinguished sequence of top elements\/} of $P$.  Clearly,
$P$ is a projective cover of any module with top $T$; in other words, the
modules with top $T$ are precisely the quotients $P/C$ with $C \subseteq
JP$, up to isomorphism.  Finally, {\it we fix a dimension vector $\bd =
(d_1,
\dots, d_n)$ of total dimension $d$\/} such that $\bt \le \bd$ in the
componentwise partial order on $\NN_0^n$.

The tops lead to a first rough subdivision of $\modlabd$:  By
$\toptbd$ we  denote the locally closed subvariety of $\modlabd$ which
consists of the points representing the modules with top
$T$, that is, 
$$\toptbd = \{x \in \modlabd \mid \text{the module corresponding to}\ x \
\text{has top}\ T\}.$$
While the varieties $\modlabd$ are connected, the
subvarieties $\toptbd$ need not be. 

As for the projective counterpart of $\modlabd$:  In \cite{\grassI}, we
defined $\grasstd$ to be the following closed subvariety
of the classical Grassmannian
$\Grfrak(\dim P - d,JP)$ of 
$(\dim P - d)$-dimensional subspaces
of the $K$-space $JP$, namely,
$$\grasstd = \{C \in \Grfrak(\dim P - d,JP) \mid C \text{\ is a\ }
\la\text{-submodule of\ } JP \}.$$ 
This variety comes with an obvious surjection
$$\grasstd \longrightarrow \{\text{isomorphism classes of\ }
d\text{-dimensional\ } \la\text{-modules with top\ } T\},$$  
sending $C$ to the class of $P/C$.  Clearly, the fibres of this map 
coincide with the orbits of the natural $\autlap$-action on
$\grasstd$.  

Here, we will consistently keep the dimension vector $\bd$ (not only the
total dimension $d$) of the considered representations fixed. The
projective counterpart to the quasi-affine variety $\toptbd$ is a closed
(and hence again projective) subvariety
$\grasstbd$ of $\grasstd$.  Namely,
$$\grasstbd = \{ C \in \grasstd \mid \underbardim P/C = \bd \}.$$
Clearly, $\grasstbd$ is stable under the $\autlap$-action, and the
corresponding orbits are in one-to-one correspondence with the
isomorphism classes of modules with top $T$ and dimension vector
$\bd$.    

We note that $\grasstbd$ is irreducible (or connected) precisely
when $\toptbd$ has this property.  This follows from the following
proposition due to Bongartz and the author, which links the ``relative
geometry'' of the $\GL_d$-stable subsets of $\toptbd$ to that of the
$\autlap$-stable subsets of $\grasstbd$.  

\proclaim{Proposition 2.1} {\rm{(}}See \cite{\GeomIV},{\rm{
Proposition C.)}} The assignment $\autlap.C \mapsto \GL_d.x$, which
pairs orbits $\autlap.C \subseteq \grasstd$ and
$\GL_d.x \subseteq \toptd$  representing the same
$\la$-module up to isomorphism, induces an inclusion-preserving bijection
$$\Phi: \{ \autlap\text{-stable subsets of\ } \grasstbd \} \rightarrow
\{\GL_d\text{-stable subsets of\ } \toptbd \}$$  
which preserves and reflects
openness, closures, connectedness, irreducibility, and types of
singularities.
\qed
\endproclaim
  
In particular, this correspondence permits transfer of information
concerning the irreducible components of any locally closed
$\GL_d$-stable subvariety of $\toptbd$ to the irreducible components of
the corresponding $\autlap$-stable subvariety of $\grasstbd$, and vice
versa.  Indeed, given that the acting groups, $\GL_d$ and
$\autlap$, are connected, all of their orbits are irreducible. 
Therefore all such irreducible components are again stable under the
respective actions, meaning that $\Phi$ restricts to a bijection between
the collection of irreducible components of $\grasstbd$ on one hand, and
that of $\toptbd$ on the other.  Similarly, information regarding
degenerations can be shifted back and forth across the bridge
$\Phi$ (see Section 4).

\subhead 2.C. The subvarieties $\laySS \subseteq \toptbd$ and $\grassSS
\subseteq \grasstbd$ \endsubhead

The radical layering of modules
provides us with a further partition of $\toptbd$ and
$\grasstbd$ into pairwise disjoint locally closed subvarieties.  (In
general, this partition fails to be a stratification in the technical
sense, however, even when $\la$ is hereditary.)  

\definition{Definition 2.2 and first comments} A {\it
semisimple sequence\/} with dimension vector $\bd$ and top
$T$ is any sequence   $\SS= (\SS_0, \dots, \SS_L)$ of semisimple modules
such that $\SS_0 = T$ and the semisimple module
$\sum_{0 \le l \le L} \SS_l$ has dimension vector $\bd$.  Moreover
(in order to rule out uninteresting sequences), we require that each $\SS_l
\le J^l P/J^{l+1} P$ for $0\le l\le L$. (This condition is intrinsic to
$\SS$, since
$P$ is the projective cover of $\SS_0$.)   Since we identify semisimple
modules with their isomorphism classes, a semisimple sequence amounts to a
matrix of discrete invariants keeping count of the multiplicities of the
simple modules in the individual slots of
$\SS$.  The radical layering of any finitely generated $\la$-module $M$
yields a semisimple sequence, namely $(M/JM, JM/J^2M,\dots,J^LM)$. As noted
in 2.B, we denote this sequence by $\SS(M)$.

Given such a semisimple sequence $\SS$, we define the following
action-stable locally closed subvarieties of
$\toptbd$ and $\grasstbd$, respectively:
$$\align
\laySS &= \{x \in \toptbd \mid x \ \text{corresponds to a module with
radical layering} \ \SS\} \\
\grassSS  &= \{ C \in \grasstbd \mid \SS(P/C) = \SS\}. \endalign$$
 Clearly, the one-to-one
correspondence $\Phi$ of Proposition 2.1, between the $\autlap$-stable
subsets of
$\grasstbd$ and the $\GL_d$-stable subsets of $\toptbd$, restricts to a
correspondence between the $\autlap$-stable subsets of $\grassSS$
and the  $\GL_d$-stable subsets of $\laySS$, which respects geometric
properties in the sense of the proposition. In particular, irreducible
components are preserved by $\Phi$ and $\Phi^{-1}$.
\enddefinition

\subhead 2.D. The projective variety $\grassbd$ 
\endsubhead  

We next embed the projective varieties $\grasstbd$ for a fixed dimension
vector $\bd$ and variable tops $T$ (with dimension vectors $\bt \le \bd$)
into a bigger variety
$\grassbd$  --  still projective  -- in order to obtain a
projective counterpart to the full classical affine variety
$\modlabd$.  The projective setting still offers advantages supplementing
those of the classical affine one, due to compactness and to the fact that
it invites methodological alternatives; indeed the different structures of
the operating groups --  reductive in the classical scenario, containing a
big unipotent radical in the Grassmannian setting  --   bring different
lines of the existing theory of algebraic group actions to bear. 
However, the asset of comparative smallness of
$\grasstbd$ and its orbits in comparison with the classical setting is
lost.  In fact, in general, we will find infinitely many copies of
$\grasstbd$ inside $\grassbd$. 

Let $Q$ be any finitely generated projective $\la$-module.  When working
with such a $Q$, in place of the projective cover $P$ of $T$, we transfer
the notation $z_r$ and $e(r)$ to $Q$. That is, we specify a distinguished
sequence of top elements $z_1,
\dots, z_m$ in $Q$, denote the primitive
idempotent norming $z_r$ by $e(r)$, and  refer to $Q =
\bigoplus_{1 \le r \le m} \la z_r$ as the distinguished projective cover
of the semisimple module $Q/JQ$.  Of particular interest will be the case
where $Q/JQ =
\bigoplus_{1 \le r \le n} S_i^{d_i}$ is the semisimple module with
dimension vector $\bd$.  In that case, we will denote the distinguished
projective cover $Q$ by $\bp$.

Any $\la$-module with top $\le Q/JQ$ is a quotient of $Q$, although $Q$
need not be its projective cover. To accomodate this fact, we simply
modify the definition of $\grasstbd$ to allow submodules of $Q$ that might
not be contained in $JQ$. 

\definition{Definition 2.3 and comments} The subset
$$\multline
[\grassbd]_Q : =  \\
 \{C \in \Grfrak(\dim Q - d,\, Q) \mid C \text{\
is a\ } \la\text{-submodule of\ } Q \ \text{with}\ \underbardim Q/C = \bd
\} \endmultline$$ 
is a closed, and hence projective, subvariety of the
classical Grassmann variety
$\Grfrak(\dim Q - d,\, Q)$.  Clearly, it parametrizes all isomorphism
classes of modules with dimension vector $\bd$ and top $\le Q/JQ$ via the
map that sends $C$ to the class of $Q/C$.  Again, we let the automorphism
group of the underlying projective module, that is,
$\Aut_\la(Q)$ in the present setting, act canonically on $[\grassbd]_Q$.
As we will see in Observation 2.4 below, the isomorphism classes of
modules parametrized here are in bijective correspondence with the
$\Aut_\la(Q)$-orbits of
$[\grassbd]_Q$.
\smallskip

In the most frequently considered special case, where $Q = \bp$, we denote
$[\grassbd]_Q$ by $\grassbd$.  So, explicitly,  
$\grassbd$ is the following projective variety:
$$\multline
\grassbd =  \\
\{C \in \Grfrak(\dim \bp - d,\, \bp) \mid C \ \text{is a}\ 
\la \text{-submodule of}\ \bp\ \text{with}\ \underbardim \bp/C = \bd
\}. \endmultline$$ 
In this latter case, the assignment $C
\mapsto \bp/C$ yields a parametrization of all modules with dimension
vector $\bd$, up to isomorphism, irrespective of their tops, and the
fibers of this parametrization coincide with the
$\bigautlap$-orbits of $\grassbd$.
\enddefinition 

\proclaim{Observation 2.4}  Let $Q$ be a finitely generated projective
$\la$-module as introduced ahead of Definition {\rm2.3}.  The isomorphism
classes of left $\la$-modules with dimension vector
$\bd$ and top $\le Q/JQ$ are in one-to-one correspondence with the
$\Aut_\la(Q)$-orbits of
$[\grassbd]_Q$, where $\Aut_\la(Q).C\subseteq [\grassbd]_Q$ corresponds to
the isomorphism class of $Q/C$.
\endproclaim

\demo{Proof}  Clearly, every module with dimension vector $\bd$ and top
contained in $Q/JQ$ is of the form
$Q/C$ for a suitable point $C \in [\grassbd]_Q$, and, given $f \in
\Aut_\la(Q)$, we have $Q/(f.C) \cong Q/C$.  

Suppose that $D$ is
a point in $[\grassbd]_Q$ such that $Q/D \cong Q/C$, via some isomorphism
$f:Q/C\rightarrow Q/D$ say.  Denote by
$\can_C$ and $\can_D$ the canonical epimorphisms $Q \rightarrow Q/C$
and $Q \rightarrow Q/D$, respectively.  Moreover, let
$\pi_C: Q'
\rightarrow Q/C$ be a projective cover of $Q/C$, and $\iota_C: Q'
\rightarrow Q$ a (necessarily split) embedding of $Q'$ into $Q$ such that
the composition $(\can_C) \circ (\iota_C) = \pi_C$.  Then $Q =
\iota_C(Q')
\oplus C'$ for a suitable direct summand $C' \subseteq C$, and $C =
\Ker(\pi_C\iota^*_C) \oplus C'$, where $\iota^*_C: \iota_C(Q')
\rightarrow Q'$ is inverse to $\iota_C$.  Since $Q/D
\cong Q/C$, there is also a projective cover $\pi_D: Q' \rightarrow Q/D$,
and symmetrically, we obtain a split embedding $\iota_D : Q' \rightarrow Q$
and a direct summand $D'$ of $D$ such that all of the above conditions
are satisfied on replacement of $C$ by $D$. In particular, $Q\cong
Q'\oplus C'\cong Q'\oplus D'$.  Consequently, $C'
\cong D'$, via some isomorphism $h':C' \rightarrow D'$ say.  

The maps $\pi_C\iota^*_C$ and $\pi_D\iota^*_D$ are projective covers of
$Q/C$ and $Q/D$, respectively, whence the isomorphism $f$ lifts to an
isomorphism $h''$ making the following diagram commute:
$$\xymatrixrowsep{3pc}\xymatrixcolsep{4pc}
\xymatrix{
\iota_C(Q') \ar[d]_{\pi_C\iota^*_C} \ar[r]^{h''} &\iota_D(Q')
\ar[d]^{\pi_D\iota^*_D} \\
Q/C \ar[r]^{f} &Q/D
}$$
Then $h''\bigl(\Ker(\pi_C\iota^*_C) \bigr) = \Ker(\pi_D\iota^*_D)$, and we
deduce that $h = h' \oplus h''$ is an automorphism of $Q$ which takes $C$
to
$D$.  In other words, $D$ lies in the $\Aut_\la(Q)$-orbit of
$C$, as required. \qed 
\enddemo

So, in particular, the $\bigautlap$-orbits of $\grassbd$ are in natural
bijective correspondence with the isomorphism classes of left
$\la$-modules having dimension vector $\bd$.
We record an analogue of Proposition 2.1 for the enlarged setting.  For
a proof, we again refer to \cite{\GeomIV},
Proposition C.
 
\proclaim{Proposition 2.5} {\rm (twin of 2.1)} The assignment $\bigautlap.C
\mapsto \GL_d.x$, which pairs orbits $\bigautlap.C \subseteq \grassbd$ and
$\GL_d.x \subseteq \modlabd$  representing the same
$\la$-module up to isomorphism, induces an inclusion-preserving bijection
$$\Phi: \{ \bigautlap\text{-stable subsets of\ } \grassbd \} \rightarrow
\{\GL_d\text{-stable subsets of\ } \modlabd \}$$  
which preserves and reflects
openness, closures, connectedness, irreducibility, and types of
singularities.  
\qed \endproclaim

\definition{Remark 2.6} In analogy with the above constructions, one may
define a projective variety $\grassd$ encoding all of the $d$-dimensional
$\la$-modules. Namely, let $Q$ be the free left $\la$-module $\la^{nd}$,
and set
$$\grassd := \{C \in \Grfrak(\dim Q - d,\, Q) \mid C \ \text{is a}\ 
\la \text{-submodule of}\ Q \}.$$
This variety is, in turn, equipped with a natural $\Aut_\la(Q)$-action,
the orbits of which parametrize the isomorphism classes of the left
$\la$-modules of dimension $d$.

As in Proposition 2.5, we then obtain a map from the action-stable subsets
of $\grassd$ to the action-stable subsets of $\modlad$, which preserves and
reflects relative geometric properties. 
Combining this observation with  the well-known fact that the connected
components of $\modlad$ coincide with the subvarieties $\modlabd$, we
obtain that the subsets
$$\{C \in \grassd \mid Q/C  \text{\
has dimension vector}\ \bold{d}\},$$ 
where $\bd$ traces the dimension vectors of total dimension $d$,
are the connected components of $\grassd$. In particular, they are
closed subvarieties of $\grassd$. 
\enddefinition

In light of Remark 2.6, it will  not impact the
usefulness of the considered varieties to restrict to
$\grassbd$ and $\modlabd$.    
    
Fix a semisimple module $T$ as specified at the
beginning of Section 2.B, and suppose $T  \le \bp/J\bp$, where, as
above, $\bp = \bigoplus_{1 \le r \le d} \la z_r$ is the distinguished
projective cover of $\bigoplus_{1\le r\le n} S_i^{d_i}$.    Then we find a
copy of
$\grasstbd$ as a (closed) subvariety of $\grassbd$:  Indeed, the
distinguished projective cover $P$ of $T$ can be identified with
$\bigoplus_{r\in\Delta} \la z_r$ for a suitable subset $\Delta\subseteq
\{1,\dots,d\}$. Then
$\grasstbd$ is evidently isomorphic to the subvariety of
$\grassbd$ consisting of the points $C$ which contain the elements $z_r$
for $r \notin \Delta$ such that, moreover, $C \cap \bigoplus_{r \in
\Delta} \la z_r \subseteq \bigoplus_{r \in \Delta} Jz_r$.  In fact, if
$T$ contains simple summands with multiplicity $\ge 2$, the
encompassing variety $\grassbd$ contains infinitely many copies
of $\grasstbd$.  This makes it far more unwieldy than the parametrizing
varieties $\grasstbd$ and $\grassSS$, and we will retreat to the latter
whenever possible.

Following the model of $\toptbd$ and $\grasstbd$, we further subdivide
the varieties $\grassbd$ in terms of tops and radical layerings. 

\definition{Definition 2.7 and Comments} Let $T$ and $\bp$ be as above, and
let $\SS = (\SS_0, \dots,
\SS_L)$ be a semisimple sequence with dimension vector $\bd$. (Remember
that by choice of $\bp$, the top $\SS_0$ of $\SS$ must be $\le \bp/J\bp$.) 
We set
$$\align
\biggrasstbd &:= \{C \in \grassbd \mid \bp/C \ \text{has top}\ T\}  \\
\biggrassSS &:= \{C \in \grassbd \mid \SS(\bp/C)= \SS\}.  \endalign$$ 
By Observation 2.9 below, $\biggrasstbd$ is locally closed in $\grassbd$
and, by Observation 2.11, so is $\biggrassSS$.
\smallskip  

The variety
$\biggrasstbd$ is the disjoint union of those $\biggrassSS$ for which
$\SS_0 = T$, and it contains the copies of $\grasstbd$ in $\grassbd$
mentioned above.  Similarly, the corresponding copies of
$\grassSS$ can all be found in $\biggrassSS$.
\enddefinition

We next compare the varieties $\grasstd$ and $\grassSS$ to
$\biggrasstd$ and $\biggrassSS$ in a couple of small examples.

\definition{Examples 2.8} (1)  Let $\la = K\Gamma/I$, where $\Gamma$ is the quiver  
\ignore $\xymatrixrowsep{1pc}\xymatrixcolsep{2pc} \xymatrix{  
1 \ar@/^0.2pc/[r] &2 \ar@/^0.2pc/[l] }$\endignore
and $I$ is generated by the paths of length $2$.  For
$\bd = (1,1)$, the large variety $\grassbd$ has two irreducible
components, each of which is isomorphic to a copy of $\PP^1$.  One of
these components consists the
$\bigautlap$-orbits of the modules $\la e_1$ and $S_1 \oplus
S_2$, the other component parametrizes the modules $\la e_2$ and $S_1
\oplus S_2$.  In either case, the projective module $\la e_i$ has an
$\bigautlap$-orbit isomorphic to $\AA^1$, and the orbit of $S_1 \oplus
S_2$ is a singleton.  Moreover, each component equals the closure of
$\biggrass^{S_i}_{\bd}$ in $\grassbd$.

On the other hand, $\grass^{S_i}_{\bd}$ consists of the
single point $J e_i$.
\smallskip   

\noindent (2) Let $\la = KQ$, where $Q$ is the generalized Kronecker quiver

$$\xymatrixrowsep{0.5pc}\xymatrixcolsep{4pc}
\xymatrix{
1 \ar[r] \ar@/^/[r]<0.5ex> \ar@/_/[r]<-0.5ex> &2
}$$
\medskip

\noindent and choose $\bd = (2,3)$, $T = S_1^2$, and $\SS = (S_1^2,
S_2^3)$.  The distinguished projective cover $\bp$ of $S_1^2 \oplus S_2^3$
is $\bp =
\bigoplus_{1 \le r \le 5} \la z_r$, where $\la z_r \cong \la e_1$ for $r
= 1,2$, and $\la z_r \cong \la e_2 = S_2$ for $r = 3,4,5$.  The
distinguished projective cover of $T$ is the submodule $P = \la z_1 \oplus
\la z_2$ of
$\bp$.  In particular, $J \bp$ and $J P$ are direct sums of copies of
$S_2$, and consequently every subspace is a $\la$-submodule.  In the small
Grassmannian setting one obtains
$$\grasstbd = \grassSS \cong \Grfrak(3, JP) \cong \Grfrak(3, K^6),$$
whence $\dim \grasstbd = 9$.
For the large setting, one computes: 
$$\grassbd = \Grfrak(6, e_2\bp) \cong \Grfrak(6, K^9)$$
has dimension $18$, and $\biggrasstbd = \biggrassSS$ is a dense open
subset.  Indeed, $\biggrasstbd = \{C \in \grassbd \mid \dim (C \cap
J\bp) = 3\}$.  Openness now follows from the fact that, for any point
$C \in \grassbd$, the intersection $C \cap J\bp$ has dimension at least
$3$, and therefore $\biggrasstbd$ coincides with the open subvariety  
$\{C \in \grassbd \mid \dim(C \cap J\bp) <
4\}$.  In particular, $\biggrasstbd$ has in turn dimension $18$.         
\qed \enddefinition

\subhead  2.E. Useful closed subsets of $\grassbd$ and $\grasstbd$
\endsubhead 

The following fact is straightforward from upper
semicontinuity of the maps $C \mapsto \dim \Hom_\la(\bp/ C, S_i)$ for $1
\le i \le n$.

\proclaim{Observation 2.9}  Given any semisimple
module $T \le
\bp/ J \bp$, the subset 
$$\{C \in \grassbd \mid T \le \operatorname{top}(\bp/ C) \}$$
of $\grassbd$ is closed. Consequently, the set $\biggrasstbd$ is a locally
closed subvariety of
$\grassbd$.
\qed \endproclaim

Next, we introduce a partial order on
the (finite) set of semisimple sequences of dimension vector $\bd$.
In analogy with Observation 2.9, it will provide us with a useful
array of closed subsets of $\grasstbd$.
 
\definition{Definition 2.10} Let $\SS$ and $\SS'$ be two semisimple
sequences with the same dimension vector.  We say that
$\SS'$ {\it dominates\/} $\SS$ and write $\SS \le \SS'$ if and only if
$\bigoplus_{l \le r} \SS_l \le \bigoplus_{l \le r}
\SS'_l$ for all $r \ge 0$. 
\enddefinition 

Roughly speaking, $\SS'$ dominates $\SS$ if and only if $\SS'$ results
from $\SS$ through a finite sequence of leftward shifts of simple
summands of $\bigoplus_{0 \le l \le L} \SS_l$ in the layering provided
by the
$\SS_l$. 

In intuitive terms, the next observation says that the simple
summands in the radical layers of the modules represented by
$\biggrassSS$ are only ``upwardly mobile" as one passes to modules in
the boundary of the closure of
$\biggrassSS$ in $\grassd$.   

\proclaim{Observation 2.11}  Suppose that
$\SS$ is a semisimple sequence with dimension vector $\bd$.  
Then the union $\bigcup_{\SS' \ge \SS,\, \underline{\dim}\SS'=\bd}
\biggrass{(\SS')}$ is closed in $\grassbd$. 
Analogously, the union $\bigcup_{\SS' \ge \SS,\, \underline{\dim}\SS'=\bd}
\lay{\SS'}$ is closed in $\modlabd$. 

In particular, $\biggrassSS$ is a locally closed subvariety of $\grassbd$.
\endproclaim

\demo{Proof} For simplicity, we posit that all semisimple sequences
mentioned in the proof have dimension vector $\bd$.

By Proposition 2.5, it suffices to prove the claim in the
Grassmannian setting.  We write $J^l \bp / J^{l+1} \bp = \bigoplus_{1
\le i \le n} S_i^{t_{li}}$ and $\SS_l =
\bigoplus_{1 \le i \le n} S_i^{s_{li}}$ for each $l \in \{0, \dots,
L\}$;  note that $s_{li} \le t_{li}$, since $\SS_l \le J^l \bp /
J^{l+1} \bp$ for all $l \ge 0$ by the definition of a semisimple
sequence with a top contained in $\bp/J\bp$.  Moreover, for each vertex
$e_i$ of $\Gamma$, we consider the following partial flag of subspaces of
the $K$-space $\bp$:
$$e_i J^L \bp \subseteq e_i J^{L - 1} \bp \subseteq \cdots
\subseteq e_i J \bp \subseteq e_i \bp.$$ 
Then the union of the $\biggrass{\SS'}$, where $\SS'$
traces the semisimple sequences dominating
$\SS$, coincides with $\grassbd \cap 
\bigl(\bigcap_{1 \le i \le n} V_i\bigr)$, where
$$V_i = \{C \in \Grfrak(\dim \bp - d,\, \bp) \mid \dim \bigl( C \cap e_i
J^l \bp \bigr)  \ge  \sum_{k \ge l} \bigl( t_{ki} - s_{ki} \bigr) \text{\
for all\ } 1 \le l \le L\}.$$  
The latter sets are well known to be closed in
$\Grfrak(\dim \bp - d,\, \bp)$, and our first claim follows.

To derive that $\biggrassSS$ is locally closed in $\grassbd$, let
$\SS^{(1)},\dots,\SS^{(u)}$ be the semisimple sequences which are strictly
larger than $\SS$. Then
$$\biggrassSS= \biggl(\; \bigcup_{\SS' \ge \SS} \biggrass{(\SS')} \biggr)
\setminus \biggl(\; \bigcup_{1\le i\le u} \ \bigcup_{\SS' \ge \SS^{(i)}}
\biggrass{(\SS')} \biggr) .$$
Since there are only finitely many semisimple sequences of any given
dimension vector, this proves our final claim.
\qed
\enddemo

\proclaim{Corollary 2.12}  {\rm{(}}cf\.
\cite{\degen}{\rm{, Observation 3.4.)}} Suppose that
$\SS$ is a semisimple sequence with dimension vector $\bd$ and top
$T$.  Then the union 
$$\bigcup_{\SS' \ge \SS, \,\SS'_0 = T, \, \underline{\dim}\, \SS' = \bd}
\grass{(\SS')}$$
 is closed in $\grasstd$. {\rm {(}}As before, $\SS'_0$ is
the semisimple module in the $0$-th slot of the sequence $\SS'$.{\rm
{)}}\qed
\endproclaim

\subhead 2.F.  Structure of the acting automorphism groups \endsubhead

The presence of a, usually large, unipotent radical in the acting
automorphism groups $\Aut_\la(Q)$ (e.g., for $Q = P$ or $Q = \bp$) will
turn out to be of advantage in making a different arsenal of methods
applicable than are available for the action of the reductive group
$\GL_d$ in the classical affine setting.  The following facts are
straightforward, but very useful.

\proclaim{Proposition 2.13 and terminology} Let $Q$ be any finitely
generated projective $\la$-module.  Then the automorphism group
$\Aut_\la(Q)$ is isomorphic to the semidirect product 
$$\bigl(\Aut_\la(Q)\bigr)_u \rtimes \Aut_\la(Q/JQ),$$
where $\bigl(\Aut_\la(Q)\bigr)_u$ is the unipotent radical of
$\Aut_\la(Q)$.

Moreover, $\bigl(\Aut_\la(Q)\bigr)_u$ equals
$\{\id_Q + g \mid g \in \Hom_\la(Q,JQ)\}$. 
\smallskip

If $Q = \bigoplus_{1 \le r \le m} \la z_r$ with distinguished top elements
$z_r$, then the following subgroup $\T_0$ of
$\Aut_\la(Q)$ {\rm (}and of $\Aut_\la(Q/JQ)$, up to isomorphism{\rm )} will
be referred to as the {\rm distinguished torus} in $\Aut_\la(Q)$:  Namely,
$$\T_0 = \{f \in \Aut_\la(Q) \mid z_1,\dots,z_m \
\text{are eigenvectors of}\ f\}. \qquad\square$$
\endproclaim

As a conseqence, we can split up the investigation of orbits and orbit
closures in the above setting into two tasks: exploring the action of the
unipotent group
$\bigl(\Aut_\la(Q)\bigr)_u$ and that of the reductive group $\Aut_\la(Q/JQ)
\cong \prod_{1
\le i \le n} \GL_{m_i}(K)$, where $m_i$ is the multiplicity of the
simple module $S_i$ in the top of $Q$.

In particular, the following results will be
useful on many occasions.  The first is due to Rosenlicht \cite{\RosII},
the second was first proved by Kostant for $K = \CC$ and subsequently
generalized to arbitrary base fields by Rosenlicht \cite{\Ros}, Theorem 2;
cf\. \cite{\Bor}, Proposition 4.10.

\proclaim{Theorem 2.14}  Let $G$ be a connected unipotent algebraic
group and
$V$ a variety on which $G$ acts morphically. Then: 

{\rm (1)} All $G$-orbits of $V$ are isomorphic to affine spaces.

{\rm (2)}  If $V$ is quasi-affine, then all $G$-orbits are closed in $V$.
\qed
\endproclaim 

To apply part (1), we typically take $V$ to be the variety $[\grassbd]_Q$,
where
$Q$ is a projective module with distinguished sequence $z_1, \dots,
z_m$ of top elements as in Definition 2.3 and Proposition 2.13, and
let $G =
\bigl(\Aut_\la(Q) \bigr)_u$ be the unipotent radical of $\Aut_\la(Q)$. 
For $C \in [\grassbd]_Q$, we spell out an explicit isomorphism of
varieties $\bigl(\Aut_\la(Q) \bigr)_u. C \rightarrow \AA^\m$, where $\m$
is the dimension of the considered orbit. 
Setting $M = Q/C$, we obtain in analogy to \cite{\degen}, Observation 3.2:
$$\m = \dim \Hom_\la (Q, JM) - \dim \Hom_\la(M, JM).$$  
The isomorphism we describe will be used in the proof of Proposition 4.4. 
Note that the stabilizer of $C$ in $\bigl(\Aut_\la(Q)
\bigr)_u$ equals the set of all automorphisms of the form $\id_Q + g$,
where $g$ runs through the subspace 
$$\stab_{\Hom_\la (Q, JQ)}(C) = \{g \in \Hom_\la(Q, JQ) \mid g(C) \subseteq
C\}.$$  
Pick homomorphisms $g_1, \dots, g_\m \in \Hom_\la (Q, JQ)$ which consitute
a $K$-basis for $\Hom_\la(Q,JQ)$ modulo $\stab_{\Hom_\la (Q, JQ)}(C)$. 
Then
$$\AA^\m \rightarrow \bigl(\Aut_\la(Q) \bigr)_u.
C, \qquad (a_i) \mapsto (\id_Q + \sum_i a_i g_i).C$$
is an isomorphism as desired.  To describe its image in concrete
computations, it is often advantageous to content oneself with a generating
set
$(g_i)$ for
$\Hom_\la(Q,JQ)$ modulo $\Hom_\la (Q, C\cap JQ)$.  Particularly convenient
are the maps $g_{r,i}$ with $g_{r,i} (z_s) =  p_{r,i} z_{j(r,i)}$ for $s =
r$ and $g_{r,i} (z_s) = 0$ for $s \ne r$, where $p_{r,i} \zhat_{j(r,i)}$
trace those paths in a skeleton of $M = Q/C$ which have positive length
and end in $e(r)$.

\head 3. Skeleta  of modules and affine charts for $\grassbd$ and
$\grasstbd$ \endhead

The affine subvarieties $\biggrassS$ covering $\grassbd$ which we will
next introduce and scrutinize are generalizations of those that were
previously defined in \cite{\generic} for $\grasstbd$; predecessors for
special choices of $T$ were considered in \cite{\GeomII, \grassI, \degen}.
These affine charts are distinguished by the following three assets:
First, they are always stable under the action of the unipotent radical
$\bigunirad$ of the acting automorphism group.  This
brings the theory of unipotent group actions to bear on the individual
patches; see, e.g., Theorem 2.14.  In fact, the
$\biggrassS$ are even stable under the action of
$\bigunirad \rtimes \T_0$, where $\T_0$ is the distinguished maximal torus
in $\bigautlap$ introduced in Proposition 2.13.  Another plus lies in the
tight connection between these affine charts and the structural features
of the modules they represent.  The $\biggrassS$
consist precisely of the points in $\grassbd$  which correspond to
modules sharing a ``special path basis" $\S$, which we will call a {\it
skeleton\/}; such a skeleton pins down other structural data, such as the
radical layering, of the pertinent module.  Finally, for large classes of
algebras beyond the hereditary ones, the $\biggrassS$ are rational.  

The affine cover
$(\biggrassS)_\S$ of $\grassbd$ will turn out to consist of
intersections of Schubert cells with the varieties
$\biggrassSS$.  In particular, we hence find it to be the union of affine 
covers of the individual subvarieties $\biggrassSS$.  (On the side, we
mention that, while the classical affine Schubert cells of
$\Grfrak(\dim \bp - d,\, \bp)$ obviously intersect $\grassbd$ in affine
sets  --  the intersection with the Schubert cell that accompanies
any skeleton $\S$ is labeled $\bigSchu(\S)$ below  --  these intersections
are hardly ever
$\bigunirad$-stable; nor are they stable under other relevant subgroups of
$\bigautlap$ in general.)  We point out that
skeleta of modules are generic attributes in the following sense:  Every
irreducible component of the variety
$\modlabd$ contains a dense open subset such that all modules in this
subset have the same skeleton (\cite{\generic}; combine Observations 2.1
and 2.4 with Theorem 3.8).  Finally, there is a fast algorithm
\cite{\codes} which provides the polynomials defining the distinguished
affine patches, starting from a presentation of
$\la$ in terms of quiver and relations.  The package is coupled with an
algorithm that decomposes the varieties $\grassS$ into irreducible
components.  It is on the level of the
$\biggrassS$ covering
$\biggrassSS$ (and the $\grassS$ covering $\grassSS$) that links
between algebraic features of the modules being
parametrized on one hand, and the geometry of the parametrizing variety
on  the other, become most accessible.  

A downside of the $\biggrassS$: In general, these affine
patches are not open in $\grassbd$, but only open in
the subvarieties $\biggrassSS$.  Moreover, the corresponding
$\bigautlap$-orbits have comparatively large dimension and are far less
manageable than the $\autlap$-orbits of the variety $\grasstbd$, where
$T$ equals the top $\SS_0$ of $\SS$.  One counters it by moving back and
forth between the large and small settings. 
  
\subhead 3.A.  Definition of skeleta and a first example \endsubhead

Roughly speaking, skeleta allow us to carry some of the benefits of the
path-length grading of the projective $K\Gamma$-modules to arbitrary
$\la$-modules.
\smallskip

{\bf Setting:} As in Section 2.B, we fix a semisimple module $T
\in \lamod$ with dimension vector $\bt$ and total dimension
$t$.  Moreover,
$\bd$ continues to denote  a dimension vector $\ge \bt$ of total dimension
$d$.  As in Section 2.D, we fix a projective
$\la$-module $Q$ such that $T \le Q/JQ \le \bigoplus_{1 \le i
\le n} S_i^{d_i}$, and pair it with a distinguished sequence of
top elements $z_1, \dots, z_m$ for $Q$ such that $z_r = e(r) z_r$ for
primitive idempotents $e(r) \in \{e_1, \dots, e_n\}$.  In the
special cases where
$Q$ is the distinguished projective cover of $T$, or of the semisimple
module
$\bigoplus_{1 \le i \le n} S_i^{d_i}$ of dimension vector $\bd$, we denote
$Q$ by $P$ or $\bp$, respectively.  

Now write $\Qhat = \bigoplus_{1 \le r \le m} K\Gamma \zhat_r$ for
the projective $K\Gamma$-module covering the projective
$\la$-module $Q =
\bigoplus_{1 \le r \le m} \la z_r$; here each
$\zhat_r$ is normed by the same vertex of the quiver $\Gamma$ as $z_r$.  In
other  words, $\Qhat / I \Qhat \cong Q$ under the canonical
$K\Gamma$-epimorphism $\Qhat \rightarrow Q$ which sends $\zhat_r$ to
$z_r$.  In the two most frequently considered instances, we write
$\Phat = \bigoplus_{1 \le r \le t} K\Gamma \zhat_r$ and 
$\bphat = \bigoplus_{1 \le r \le d} K\Gamma \zhat_r$, respectively.  

By a {\it path in $\Qhat$ starting in $\zhat_r$\/} we mean any element $p'
= p \zhat_r \in \Qhat$, where $p$ is a path in $K\Gamma$ starting in the
vertex $e(r)$.  The {\it length\/} of $p'$ is defined to be that of $p$;
ditto for the {\it endpoint\/} of $p'$.  If $p_1$ is an initial subpath of
$p$, meaning that $p = p_2 p_1$ for paths $p_1$, $p_2$, we call
$p_1\zhat_r$ an {\it initial subpath of\/} $p \zhat_r$.  So, in
particular, $\zhat_r = e(r)\zhat_r$ is an initial subpath of length $0$
of any path $p \zhat_r$ in
$\Qhat$. 
\smallskip 

The reason why we do not identify $p' = p \zhat_r$ with
$p$ lies in the fact that we want to distinguish between $p \zhat _r$
and $p \zhat_s$ for $r \ne s$ but $e(r) = e(s)$.  A priori, we moreover
distinguish between a path $p \zhat_r \in \Qhat$ and the
corresponding residue class $p z_r \in Q = \Qhat/ I \Qhat$ in order to
keep an unambiguous notion of length.  However, in the sequel, we will
often not uphold the distinction between
$p \zhat_r$ and $p z_r$, unless there is a need to emphasize
well-definedness of path lengths. 

We always have distinguished embeddings $P \subseteq Q \subseteq \bp$
(those respecting the distinguished top elements carrying the same label),
and correspondingly, $\Phat \subseteq \Qhat \subseteq \bphat$. Some of the
theory we explicitly develop only for the special cases
$Q = P$ and $Q = \bp$, in order to pare down the setting to a more tightly
determined one for clarity.  Unless we emphasize the contrary, the results
carry over to the general situation with analogous proofs.   

\definition{Definition 3.1 and Conventions}
  
\noindent {\bf {(1)}} An ({\it abstract\/}) {\it skeleton with  top $T$} in
$\Qhat$ is a set $\S$ of paths of lengths at most $L$ in
$\Qhat$, which 
contains $t$ paths of length zero, say $\zhat_{i_1}, \dots,
\zhat_{i_t}$ such that $\bigoplus_{1 \le r
\le t} \la z_{i_r}$ is a $\la$-projective cover of
$T$; moreover, $\S$ is required to be closed under initial subpaths in the
following sense: whenever
$p_2 p_1 \zhat_r \in \S$, then $p_1 \zhat_r \in \S$.
(We usually suppress reference to the $K\Gamma$-projective module
$\Qhat$, when our choice of the $\la$-projective module
$Q$ is clear.)  

We say that $\S$ has {\it dimension vector\/}
$\bd$ if, for each $i$ between $1$ and $n$, the set
$\S$ contains exactly $d_i$ distinct paths ending in the vertex $e_i$.   
\smallskip   

\noindent Alternatively, we view a skeleton
$\S$ as a forest of $t$ tree graphs, where each tree consists of the paths
in
$\S$ starting in a fixed top element $\zhat_r$; see Example 3.4 below and
the remarks preceding it.  

Note:  In case $Q = P$ is the distinguished projective cover of $T$,
any abstract skeleton with top $T$ in $\Phat$ contains the full
collection of paths of length zero in $\Phat$. In light of the canonical
embeddings 
$\Phat \subseteq \Qhat \subseteq \bphat$, we routinely view
a skeleton in $\Phat$  as a subset of $\bphat$, whenever convenient.   
\smallskip

Further conventions:  
 
\noindent $\bullet$  By $\S_l$ we denote the set of paths of length $l$
in $\S$.

\noindent $\bullet$  When we pass from the $K\Gamma$-module $\Qhat$
to the $\la$-module $Q$ by factoring out $I \Qhat$, we
identify $\zhat_{i_r}$ with $z_{i_r}$ and view $\S$ as the subset 
$$\{p z_{i_r} \in Q \mid  p \zhat_{i_r} \in \S,\, r \le t \}$$ 
of $Q$ 
\smallskip    

\noindent {\bf {(2)}}  Let $\S$ be an abstract skeleton
with top $T$ in $\Phat$ (or, more generally, in $\Qhat$), and $M$ a
$\la$-module.  

We call $\S$ a {\it skeleton
of\/} $M$ if there exists a full sequence of top elements $x_1, \dots,
x_t$ of $M$, together with a $K\Gamma$-epimorphism $f: \Phat \rightarrow
M$ (resp.,
$\Qhat \rightarrow M$) satisfying $f(\zhat_{i_r}) = x_r$ for all $r$ such
that, for each $l \in \{0, \dots, L\}$, the set
$$f(\S_l) = \{f(p \zhat_{i_r}) \mid p \zhat_{i_r} \in \S_l\}  =  \{(p
+ I) x_r \mid r \le t, \, p \zhat_{i_r} \in \S_l\}$$  
induces a $K$-basis for the radical layer $J^lM / J^{l+1} M$.  In this
situation, we also say that $\S$ is a {\it skeleton of $M$ relative to
the sequence $x_1, \dots, x_t$ of top elements\/}, and observe that the
union over $l \le L$ of the above sets $f(\S_l)$, namely $\{(q + I)x_r
\mid  r \le t, \, q \zhat_{i_r} \in \S \}$, is a $K$-basis for
$M$.  Thus the skeleta of $M$ are $K$-bases which
are closely tied to the quiver presentation of $\la$.  

If $M$ has skeleton $\sigma \subseteq \Qhat$, then clearly $M$ is an
epimorphic image of $Q$.  If $M = Q/C$, we call
$\S$ a {\it distinguished skeleton of
$M$\/} provided that
$\S$ is a skeleton of $M$ relative to the distinguished sequence $z_{i_1} +
C, \dots, z_{i_t} + C$ of top elements of $Q/C$.
\smallskip    

\noindent {\bf {(3)}} If $Q = P$ is the projective cover of $T$, we write
$$\grassS = \{ C \in \grasstd \mid \S \text{\ is a distinguished
skeleton of\ } P/ C \},$$
as in \cite{\generic}.

If, on the other hand, $Q = \bp$ is the  projective cover of
$\bigoplus_{1 \le i \le n} S_i^{d_i}$, we define 
$$\biggrassS = \{ C \in \grassbd \mid \S \text{\ is a distinguished
skeleton of\ }  \bp/C \}.$$ 
\enddefinition

As we pointed out before, $\grassS$ and $\biggrassS$ parametrize the
same collection of isomorphism types of modules.  The set of all
skeleta of a $d$-dimensional module
$M$ is clearly an isomorphism invariant of $M$.  This set is always
nonempty and finite, the latter being due to our requirement that the
distinguished projective modules
$Q$ have tops contained in the semisimple module $(\la/J)^d$.  Moreover, if
$M = Q/C$, this set contains at least one distinguished candidate.  We
record a slightly upgraded variant of the existence statement, leaving
the easy proof to the reader.

\proclaim{Observation 3.2}  Suppose $M \in \lamod$ has top $T$, and let
$x_1, \dots, x_t$ be a full sequence of top elements of $M$.  Then $M$
has at least one skeleton in $\Phat$ relative to $x_1, \dots, x_t$. 
Moreover, whenever a module $\tilde{M}$ has the same
top as $M$ and maps epimorphically onto
$M$, any skeleton of
$M$ can be extended to a skeleton of $\tilde{M}$.  In
particular, every skeleton of $M$ is contained in a skeleton of $P$. 

Analogously, any module $M$ with dimension vector $\bd$ has at least
one skeleton in $\bphat$, and every such skeleton is contained in a
skeleton of $\bp$ in $\bphat$. \qed
\endproclaim  

Our definition of a skeleton coincides in
essence with that given in \cite{\grassI} and \cite{\degen} for the
situation $Q = P$ and a squarefree top
$T = P/JP$.  However, in that special case, it is unnecessary to hook up
the elements of an abstract skeleton $\S$ with top elements of
the $K\Gamma$-module $\Phat$, since the dependence on specific sequences of
top elements disappears.  In particular,
every skeleton of $P/C$ is distinguished in the case of squarefree
$T$:  Indeed, if $y_1, \dots, y_t$ is any sequence of top elements
of $P$ such that each $y_i$ is normed by the same primitive
idempotent as $z_i$, the $y_i$ can differ from the $z_i$ only by a
nonzero constant factor and a summand in $JP$; consequently, any
skeleton of $P/C$ relative to $y_1, \dots, y_t$ is also a skeleton
relative to the distinguished sequence $z_1, \dots, z_t$.        

Again, let $Q$ be either $P$ or $\bp$, and let $\S$ be an abstract skeleton
with dimension vector $\bd$ and an arbitrary top $T$.  Suppose that
$M = Q/C$.  Then $\S$ is a (not
necessarily distinguished) skeleton of $M$ if and only if $M \cong Q/D$ for
{\it some\/} point $D$ in $\grassS$ (resp., in $\biggrassS$).  However,
$\autlap.C \cap
\grassS \ne \varnothing$ need not imply $C \in \grassS$, and an analogous
caveat pertains to the big setting.  On the other hand, the
$\grassS$ (resp., $\biggrassS$) do enjoy the following partial stability
under the $\autlap$-action (resp., $\bigautlap$-action):

\definition{Observation 3.3 and Convention} The 
subvarieties $\grassS$ of $\grasstbd$ are stable under the action
of $G = \unirad \rtimes \T_0$, where $\unirad$ is the unipotent radical of
$\autlap$ and $\T_0$ the distinguished torus in
$\autlap$ (see Proposition 2.13); we write the elements of $\T_0$ in the
form $(a_1,\dots, a_t) \in (K^*)^t$, the latter standing for the
automorphism that sends $z_r$ to $a_r z_r$. 

Analogously, the sets $\biggrassS$ are  stable under the action
of the group $G = \bigunirad \rtimes {\bold \T}_0$, where $\bigunirad$ and
${\bold
\T}_0$ are the unipotent radical and distinguished torus in $\bigautlap$.
\enddefinition

\demo{Proof}  We only address the small scenario $\grasstbd$, the
arguments for the big being analogous.  Suppose $C \in \grassS$, meaning
that $\S$ is a distinguished skeleton of $P/C$.  Then $\S$ remains a
skeleton of $P/C$ after passage from our given
distinguished sequence $z_1, \dots, z_t$ of top elements of $P$ to a
new sequence of the form $g.z_1, \dots, g.z_t$ for any
$g \in G$.  But this is tantamount to saying that $\S$ is a
distinguished skeleton of $P/(g^{-1}.C)$ relative to the original
sequence.  In other words, $g^{-1}.C \in \grassS$.
\qed \enddemo  

Any abstract skeleton $\S$ can be communicated by means of an undirected
graph which is a {\it forest\/}, that is, a finite disjoint union of
tree graphs:  There are $t$ trees if $\S$ is a skeleton with top $T$,
one for each $r \le t$;  here $\zhat_r$, identified with $e(r)$,
represents the root of the $r$-th tree, recorded in the top row of the
graph. The paths $p \zhat_r$ of positive length in $\S$ are represented
by edge paths of positive length.  Instead of formalizing this
convention, we will illustrate it in Example 3.4(2).  

In general, neither $\grasstbd$ nor $\grassbd$
has an affine cover which is closed under the full $\autlap$-action (resp.,
$\bigautlap$-action), as witnessed by the first of the following examples. 
Indeed, the existence of such a cover would force all orbits to be
quasi-affine.  

\definition{Examples 3.4}  

\noindent (1) Let $\Gamma$ be the quiver $1 \rightarrow 2$, $\bd =
(2,1)$, $T = S_1^2$ and $\SS = (S_1^2, S_2)$.  Then
$\grasstbd \cong \PP^1$ consists of a single $\autlap$-orbit, where $P =
\la z_1 \oplus \la z_2 \cong (\la e_1)^2$.  Here $\bp = P
\oplus \la z_3$, where $\la z_3 = \la e_2 = S_2$, and $\grassbd \cong
\PP^2$ consists of two $\bigautlap$-orbits.  One of them equals
$\biggrassSS$, the other is a singleton, namely $\{J \bp\}$.  Clearly,
the big orbit again fails to be quasi-affine. 

\noindent (2)  In the following, we give an example of a module with
several distinct skeleta, two of which we will present in graphical
format.  Let
$\la = K\Gamma/I$, where $\Gamma$ is the quiver with a single vertex,
labeled $1$, and three loops, $\alpha$, $\beta$, $\gamma$, and $I \subseteq
K\Gamma$ the ideal generated by $\alpha^2$, $\beta^2$, $\gamma^2$, and all
paths of length $4$.

We denote by $P = \bigoplus_{1 \le r \le 3} \la z_r$ the distinguished
projective cover of $T = S_1^3$, let $d = 11$, and consider the
$d$-dimensional module $M = P/C$, where $C \in \grasstd$ is generated by
$\gamma z_1$, $\gamma \alpha z_1$, $\beta z_2$, $\gamma z_2$, $\gamma
\alpha z_2$, $\alpha z_3$, $\beta z_3$, $(\beta \alpha - \alpha \beta)z_1$,
$\alpha \gamma \beta z_1 - \beta \alpha z_2 - \gamma z_3$, and all
elements of the form $p z_2$, $q z_3$, where $p$ is a path of length $3$
and $q$ a path of length $2$ in $\Gamma$.  The module $M$ can be visualized
as follows:

$$\xymatrixrowsep{1.5pc}\xymatrixcolsep{1pc}
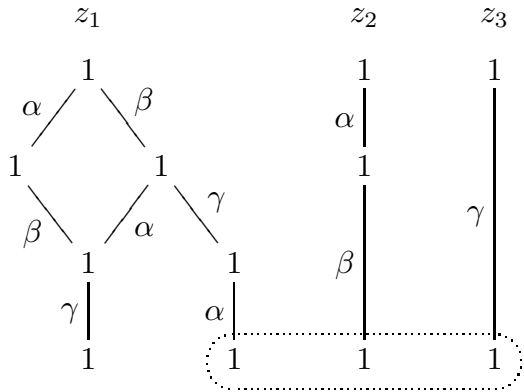
\xymatrix{
 &1 \dropup{z_1} \edge[dl]_{\alpha} \edge[dr]^{\beta} &&&&1 \dropup{z_2}
\edge[d]_{\alpha} &&1 \dropup{z_3} \edge[ddd]_{\gamma} \\
1\edge[dr]_{\beta} &&1 \edge[dl]^{\alpha} \edge[dr]^{\gamma} &&&1
\edge[dd]_{\beta} \\
 &1 \edge[d]_{\gamma} &&1 \edge[d]_{\alpha} \\
 &1 &&1 \levelpool{4} &&1 &&1 \dropdown{}
}$$

\noindent Here we use the graphing technique introduced formally in
\cite{\generic} and informally in \cite{\menace}. In particular, the
grouping of vertices inside the dotted curve means that the images in $M$
of the paths
$\alpha\gamma\beta z_1$, $\beta\alpha z_2$, $\gamma z_3$ in $P$ are
$K$-linearly dependent, whereas any two of them are $K$-linearly
independent.   Up  to permutations of the top elements
$\zhat_r$ of the projective
$K\Gamma$-module $\Phat =
\bigoplus_{1 \le r \le 3} K\Gamma \zhat_r$, the module
$M$ has precisely two skeleta in $\Phat$.  Namely,

$$\xymatrixrowsep{1.5pc}\xymatrixcolsep{0.5pc}
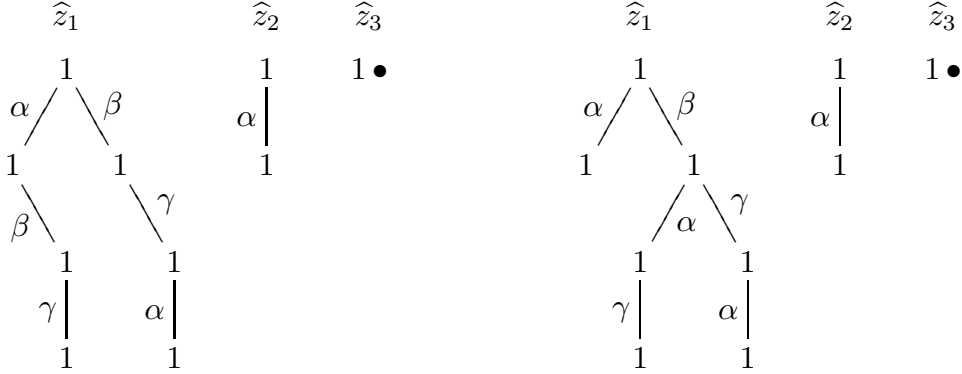
\xymatrix{
 &1 \dropup{\zhat_1} \edge[dl]_{\alpha} \edge[dr]^{\beta} &&&&1 \dropup{\zhat_2}
\edge[d]_{\alpha} &&1\,\bullet \dropup{\zhat_3} 
&&&&& &1 \dropup{\zhat_1} \edge[dl]_{\alpha} \edge[dr]^{\beta} &&&&1 \dropup{\zhat_2}
\edge[d]_{\alpha} &&1\,\bullet \dropup{\zhat_3} \\
1 \edge[dr]_{\beta} &&1 \edge[dr]^{\gamma} &&&1 
&&&&&&&1 &&1 \edge[dl]^{\alpha} \edge[dr]^{\gamma} &&&1 \\
 &1 \edge[d]_{\gamma} &&1 \edge[d]_{\alpha} 
&&&&&&&&& &1 \edge[d]_{\gamma} &&1 \edge[d]_{\alpha} \\
 &1 &&1  &&&&&&&&& &1 &&1 
}$$

\noindent Each of the two skeleta consists of all the paths $p\zhat_r$ that
occur as edgepaths in one of the $3$ trees on reading them from top to
bottom.  \qed \enddefinition

\subhead 3.B.  $\grassSS$ and $\biggrassSS$ as unions of
$\grassS$'s and $\biggrassS$'s, respectively \endsubhead 

Suppose that $\SS$ is a semisimple sequence with dimension vector $\bd$ and
top $T$.  We continue to denote by $P$ and $\bp$ the distinguished
projective covers of $T$ and $\bigoplus_{1 \le i \le n} S_i^{d_i}$,
respectively.  We next give a necessary condition for nontriviality of the
intersections $\grassS \cap \grassSS$ and $\biggrassS \cap \biggrassSS$.  

\definition{Definition 3.5} Keep the notation of Definition 3.1.  Given a
semisimple sequence $\SS$ =  $(\SS_0, \dots, \SS_L)$, we call a skeleton
$\S \subseteq \Qhat$ {\it compatible with $\SS$\/} if, for each
$l \le L$ and $i \le n$, the number of paths in $\S_l$ ending in the
vertex $e_i$ coincides with the multiplicity of the simple module $S_i$
in $\SS_l$. 
\enddefinition

In case $Q = P$, this compatibility concept
coincides with the one introduced in \cite{\generic}.
If $\grassSS \ne \varnothing$ (or, equivalently, $\biggrassSS \ne
\varnothing$), there generally are numerous (even though only finitely
many) skeleta compatible with $\SS$.  On the other hand, there is only one
semisimple sequence $\SS$ with which a given skeleton is compatible.

Clearly, each abstract skeleton in $\Qhat$ which is compatible with
$\SS$ shares top and dimension vector with $\SS$.  Moreover, given a
$\la$-module $M$, each skeleton of
$M$ is compatible with the radical layering $\SS(M)$ of $M$; this is an
immediate consequence of the definitions.  Consequently, compatibility of
$\S$ with
$\SS$ is a necessary condition for the intersection $\grassS \cap \grassSS$
(resp.,
$\biggrassS \cap \biggrassSS$) to be nonempty.  More precisely:

\definition{Observation 3.6}  The following statements are equivalent for
$M \in \lamod$ and a skeleton $\S$ in $\Qhat$:

\noindent $\bullet$ $\S$ is a skeleton of $M$.

\noindent $\bullet$ $\S$ is compatible with
$\SS(M)$, and there exists a $K\Gamma$-epimorphism $f: \Qhat
\rightarrow M$ such that $f(\S)$ is a $K$-basis for $M$. \qed  
\enddefinition

In particular, we glean:  If $\grassS \ne \varnothing$, then $\S$ is
compatible with $\SS$ precisely when  
$\grassS \subseteq \grassSS$ (equivalently, when $\biggrassS \subseteq
\biggrassSS$).

If $Q = P$, we denote by
$\Schu(\S)$ the intersection of $\grasstbd$ with the big
Schubert cell in the classical Grassmannian
$\Grfrak(\dim P - d,JP)$ that corresponds to $\S$.  In other words, 
$\Schu(\S)$ consists of all points $C \in \grasstbd$ with the property
that 
$$JP \quad = \quad C \oplus \bigoplus_{p \zhat_r \in \S,\ p\zhat_r\,
\text{of positive length}} K pz_r.$$
The big counterpart to $\Schu(\S)$ obtained for $Q = \bp$ will be denoted
by $\bigSchu(\S)$.  It lives in $\Grfrak(\dim \bp - d,\, \bp)$.  Namely,
$$\bigSchuS = \{C \in \grassbd \mid  \bp = C \oplus \bigoplus_{p \zhat_r
\in \S} Kp z_r \}.$$
In particular, $\Schu(\S)$ and $\bigSchu(\S)$ are open subsets of
$\grasstbd$ and $\biggrasstbd$, respectively.

\definition{Observation 3.7}  Suppose that $\S$ is a skeleton which is
compatible with the semisimple sequence  $\SS$.  Then 
$$\grassS =  \grassSS \cap \Schu(\S) \quad \text{and} \quad  \biggrassS = 
\biggrassSS \cap \bigSchu(\S). \ \ \qed$$
 \enddefinition   

Combining Observations 3.6 and 3.7, we obtain:

\proclaim{Corollary 3.8}  The $\grassS$ {\rm {(}}resp., $\biggrassS${\rm
{)}}, where $\S$ runs through the skeleta in $\Phat$ {\rm {(}}resp., in
$\bphat${\rm {)}} which are compatible with $\SS$, form an open cover of
$\grassSS$ {\rm {(}}resp., of $\biggrassSS${\rm {)}}.  In general, they
fail to be open in $\grasstd$ {\rm {(}}resp., $\biggrasstbd${\rm {)}},
however.  \qed
\endproclaim

To describe affine coordinates for the varieties $\grassS$ and
$\biggrassS$, we separate the small and large scenarios.

\subhead 3.C.  Critical paths and affine coordinates for the subvarieties
$\grassS$ of $\grasstd$
\endsubhead

We specialize to the situation where $Q = P$ is the distinguished
projective cover of $T$, keeping the notation of the two preceding
sections.  Moreover, we fix a skeleton $\S \subseteq \Phat$ with dimension
vector $\bd$ and top $T$.  Recall that our purpose in shifting from the
$\la$-projective module $P$ to the $K\Gamma$-projective module $\Phat$ in
defining a skeleton $\S$ was to render the lengths of the paths in
$\S$ unambiguous.  In the following supplement to Definition 3.1,
this well-definedness of path lengths is once more essential. 

\definition{Definition 3.9} A {\it $\S$-critical path\/}
is a path of length at most $L$ in $\Phat \setminus \S$, with the
property that every proper initial subpath belongs to $\S$.  In other
words, a path in $\Phat$ is $\S$-critical if and only if it fails to
belong to $\S$ and is of the form $\alpha p \zhat_r$, where 
$\alpha$ is an arrow and $p \zhat_r \in \S$.  Moreover, for any such
$\S$-critical path
$\alpha p \zhat_r$, its
$\S$-{\it set\/}, denoted $\S(\alpha p
\zhat_r)$, is defined to be the set of all paths
$q \zhat_s \in \S$ which are at least as long as $\alpha p \zhat_r$ and
end in the same vertex as $\alpha p \zhat_r$.  

Finally, we let $N$ be the {\it disjoint\/} union of the sets $\{\alpha p
\zhat_r \} \times \S(\alpha p \zhat_r )$, where $\alpha p \zhat_r$ traces
the $\S$-critical paths.  (We write the elements of $N$ as pairs, since a
priori, the sets $\S(\alpha p \zhat_r) \subseteq \Phat$ need not be
disjoint.) \qed
\enddefinition

In particular, this definition entails:  Whenever the length of a
$\S$-critical path $\alpha p \zhat_r$  exceeds the maximum of the
lengths of the paths in $\S$, the corresponding $\S$-set
$\S(\alpha p \zhat_r)$ is empty.  In accordance with our earlier
conventions (Definition 3.1), we often identify a $\S$-critical path
$\alpha p \zhat_r$ with its $I \Phat$-residue class $\alpha p z_r$ in
$P$, unless we wish to emphasize well-definedness of path lengths. If
$\S$ is the first skeleton in Example 3.4(2), we find:  The $\S$-critical
paths are $\gamma \zhat_1$, $\gamma \alpha \zhat_1$, $\alpha \beta
\alpha \zhat_1$, $\alpha \beta
\zhat_1$, $\beta \gamma \beta \zhat_1$, $\beta \zhat_2$, $\gamma \zhat_2$,
$\beta\alpha \zhat_2$, $\gamma \alpha \zhat_2$, $\alpha \zhat_3$, $\beta
\zhat_3$, $\gamma \zhat_3$.  The
$\S$-set of $\gamma \alpha \zhat_1$, for instance, is $\S(\gamma
\alpha \zhat_1)$ $=$ $\{\beta \alpha \zhat_1, \gamma \beta \zhat_1, \gamma
\beta \alpha \zhat_1, \alpha \gamma \beta \zhat_1\}$.

First, we note that the cardinality of the set $N$ does not actually
depend on $\S$. The proof is left to the reader. 

\proclaim{Observation 3.10}  Let $\S$ be a
skeleton with dimension vector $\bd$ and top $T$.  Moreover, let $\SS$ be
the unique semisimple sequence with which $\S$ is compatible. Then the
number
$$|N| = \sum_{\alpha p \zhat_r \ \sigma\text{-critical}}
|\sigma(\alpha p \zhat_r)|$$
depends only on $\SS$, not on $\S$. \qed
\endproclaim

The next observation is obvious.  We state it in order to emphasize its
pivotal role.

\definition{3.11.  Key consequence of the definition} Let $C$ be a point
in the subvariety $\grassS$ of $\grasstbd$.  Given
any $\S$-critical path $\alpha p \zhat_r$, there exist unique scalars
$c_{\alpha p \zhat_r, q \zhat_s}$ with the property that
$$\alpha p z_r + C =  \sum_{q \zhat_s \in \S(\alpha p \zhat_r)} c_{\alpha
p \zhat_r, q \zhat_s} q z_s + C$$  
in $P/C$.  The isomorphism type of $M = P/C$ is
completely determined by the resulting family of such scalars, as
$\alpha p \zhat_r$ traces all $\S$-critical paths.

Thus we obtain a map 
$$\psi: \grassS \rightarrow \AA^N,\ \ \ \ C \mapsto c =
\bigl(c_{\alpha p \zhat_r, q \zhat_s}\bigr)_{\alpha p \zhat_r\
\S\text{-critical},\ q \zhat_s \in \S(\alpha p \zhat_r)}. \ \ \qed$$ 
\enddefinition   

The principal aim of this section is to show that this map $\psi$ is an
isomorphism from the variety $\grassS$ onto a closed subvariety of
$\AA^N$.  If we identify the paths $p \zhat_r \in \Phat$ with the
corresponding elements $p z_r \in P$, we readily
find:  The point
$C \in \grassS$ corresponding to a point $c$ in the image
of $\psi$ is the $\la$-submodule of $JP$ which is generated by the elements
$\alpha p  z_r -  \sum_{q \zhat_s \in
\S(\alpha p \zhat_r)} c_{\alpha p \zhat_r, q \zhat_s} q z_s$, where
$\alpha p \zhat_r$ runs through the
$\S$-critical paths.   

The following theorem generalizes the corresponding result 
(\cite{\grassI}, Theorem 3.14) for $\grasstd$, where $T$ is a {\it
squarefree\/} semisimple module.  It serves as a fundamental tool in
\cite{\generic} and \cite{\degenII}. 

\proclaim{Theorem 3.12}  Let $\S \subseteq \Phat$ be an abstract
$d$-dimensional skeleton with top $T$. 
\smallskip

The image of the map $\psi: \grassS \rightarrow
\AA^N$ of {\rm Observation 3.11} is a closed subvariety of $\AA^N$; in
particular, it is affine.  

Moreover, $\psi$ is an isomorphism $\grassS
\rightarrow
\Im(\psi)$, whose inverse  sends any point $c =
\bigl(c_{\alpha p \zhat_r, q \zhat_s}\bigr)_{\alpha p \zhat_r\
\S\text{-critical},\ q \zhat_s \in \S(\alpha p \zhat_r)}$ to the
submodule
$U(c) \subseteq JP$ which is generated by the differences
$$\alpha p z_r -  \sum_{q \zhat_s \in \S(\alpha p \zhat_r)} c_{\alpha p
\zhat_r, q \zhat_s} q z_s,$$  
where $\alpha p \zhat_r$ runs through the
$\S$-critical paths.  This yields the following commutative diagram:
$$\xymatrixrowsep{2pc}\xymatrixcolsep{4pc}
\xymatrix{
\Im(\psi) \save+<-6.6ex,0.2ex> \drop{\AA^N\supseteq} \restore
\ar[dd]_{\psi^{-1}} \ar[dr]^-{\chi}\\
 &\save+<23ex,0ex> \drop{\txt{ \{\rm isom.~types of
$\la$-modules with skeleton $\S$\} }} \restore\\
\grassS \save+<-10.1ex,0.2ex> \drop{\grasstd\supseteq} \restore
\ar[ur]_-{\phi} }$$
 
\noindent  Here $\chi$ is the map which sends any point $c
\in\Im(\psi)$ to the isomorphism class of the factor module $P/U(c)$, and
$\phi$ the map which sends any point
$C \in \grassS$ to the class of $P/C$.
\smallskip

\noindent Polynomial equations determining the incarnation $\psi(\grassS)$
of $\grassS$ in $\AA^N$ can be algorithmically obtained from
$\Gamma$ and any set of left ideal generators for $I \subseteq
K\Gamma$.  If $K_0$ is a subfield of
$K$ over which left ideal generators for $I$ are defined, the resulting
polynomials determining $\grassS$ are defined over $K_0$ as well.
\endproclaim

Once the theorem is proved, we will identify  points $C \in \grassS$
with the corresponding points $c = \psi(C) \in \AA^N$, whenever convenient.

The algorithmic portion of the theorem will only be implicitly addressed
here. The algorithm has been implemented in \cite{\codes} with
reference to the present paper.

For a proof of Theorem 3.12, we further enlarge $K\Gamma$,
namely we consider the following (noncommutative) polynomial ring over the
path algebra $K\Gamma$, in  which the variables commute with the
coefficients from $K\Gamma$: 
$$K\Gamma[X]  = K\Gamma[X_\nu \mid \nu \in N],$$
where again
$$N = \bigsqcup_{\alpha p \zhat_r\ \S\text{-critical}} \{\alpha p
\zhat_r\} \times \S(\alpha p \zhat_r).$$  
Moreover, we extend $\Phat$ to a projective module $\P$ over the enlarged
base ring,
$$\P = K\Gamma[X] \otimes_{K\Gamma} \Phat = K\Gamma[X] \otimes_{K\Gamma}
\biggl(\bigoplus_{1 \le r \le t} K\Gamma
\, \zhat_r
\biggr) \cong \bigoplus_{1 \le r \le t} K\Gamma[X] \, \zhat_r,$$
and consider the submodule $\C$ of $\P$ which is generated over
$K\Gamma[X]$ by all differences 
$$\alpha p \zhat_r - \sum_{q \zhat_s \in \S(\alpha p \zhat_r)}
X_{\alpha p \zhat_r, q \zhat_s}\, q \zhat_s,$$
where $\alpha p \zhat_r$ runs through the $\S$-critical paths.
In particular, $\alpha p \zhat_r \in \C$ whenever $\S(\alpha p \zhat_r)$
is empty.  Roughly, the idea is to find polynomial equations for
the set of points
$c = (c_\nu)_{\nu
\in N}$ with the property that replacement of $X_\nu$ by $c_\nu$ turns
$\P/ \C$ into a left $\la$-module with $K$-basis $\S$.  (Note that $\Phat$
canonically embeds into
$\P$, whence we may view $\S$ as a subset of $\P$.) It will turn out
that we will, in fact, be solving a more general universal
problem; see Proposition 3.15.     

It will be convenient to write $y_1 \seq y_2$ whenever $y_1, y_2 \in \P$
have the same residue class modulo $\C$. 

\proclaim{Lemma 3.13}  The factor module $\P / \C$ is a free module over
the commutative polynomial ring $K[X]$, having as
basis the $\C$-residue classes of the elements in $\S$. 

In other words, given any element $y \in \P$, there exist
unique polynomials $\tau^y_{q \zhat_s} \in K[X]$ such that 
$$ y \seq \sum_{q \zhat_s \in \S} \tau^y_{q \zhat_s} \,q\, \zhat_s.$$
\endproclaim

\demo{Proof}  To show that, modulo $\C$, every element $y \in \P$ is a
$K[X]$-linear combination of elements in $\S$, it is clearly harmless to
assume that
$y$ has the form $u \zhat_r$ for some path $u \in K\Gamma$; indeed, these
elements generate $\P$ over $K[X]$.  Let $u'$ be the longest initial
subpath of
$u$ such that $u' \zhat_r \in \S$; this path $u'$ may have length
zero. If $u' = u$, then $y \in \S$. 
Otherwise, $u'$ is a proper initial subpath of $u$, and consequently
there exists a $\S$-critical initial subpath of $u
\zhat_r$, say $\alpha u' \zhat_r$, where $y = u'' \alpha u'
\zhat_r$ for a suitable terminal subpath $u''$ of $u$.  Thus $y \seq
\sum_{q \zhat_s \in \S(\alpha u' \zhat_r)} X_{\alpha u' \zhat_r, q \zhat_s}
u'' q \zhat_s$.  Each of the paths $u'' q \zhat_s$ occurring in this sum
has an initial subpath in $\S$ which is strictly longer than $\len(u')$. 
Repeat the procedure with these paths.  Since every $\S$-critical path
$\beta q \zhat_s$ of length
$\ge L+1$ has empty $\S$-set $\S(\beta q \zhat_s)$, every path whose
length exceeds $L$ is congruent to zero.  Therefore, the procedure
terminates in a $K[X]$-linear combination of paths in $\S$.

To prove linear independence over $K[X]$ of the residue classes $p
\zhat_r + \C \in \P/\C$ with $p \zhat_r \in \S$, we assume the contrary.
This amounts to the existence of distinct paths $p_1 \zhat_{r_1},
\dots, p_m \zhat_{r_m}$ in $\S$, together with nonzero elements $a_i \in
K[X]$, such that the sum $\sum_{1 \le i \le m} a_i p_i \zhat_{r_i}$
belongs to $\C$. That the listed paths are distinct in $\P$ means that,
for $i \ne j$, either $p_i \ne p_j$ or $r_i \ne r_j$.  In light of the
definition of $\C$, our assumption translates into an equality in $\P$ as
follows:
$$\multline \sum_{1 \le i \le m} a_i\, p_i\, \zhat_{r_i}  \ \ = \\ 
\sum_{ \alpha p \zhat_r \ \S\text{-critical}} \   \sum_{1 \le j \le
m(\alpha p \zhat_r)}
\ b_{\alpha p \zhat_r, j}\, u_{\alpha p \zhat_r, j}
\biggl(\alpha p \zhat_r \ -
\sum \Sb q \zhat_s\in \S(\alpha p \zhat_r) \endSb X_{\alpha
p \zhat_r, q \zhat_s}\, q \zhat_s 
\biggr), \endmultline \tag\dagger$$  
for suitable elements
$b_{\alpha p \zhat_r, j} \in K[X]$, distinct paths $u_{\alpha p
\zhat_r, j}$ in $\Gamma$ of lengths $\ge 0$ starting in the endpoint of
$\alpha$, and nonnegative integers $m(\alpha p \zhat_r)$.  Since
$\P$ is a free $K[X]$-module having as basis all paths $v \zhat_r$ in the
projective $K\Gamma$-module
$\Phat$, we are now in a position to compare coefficients.  Pick a
$\sigma$-critical path $\beta v \zhat_\mu \in \Phat$ of minimal length
such that some $b_{\beta v \zhat_\mu, j_0} \ne 0$; here $\beta$ is an arrow
and
$v \zhat_\mu \in \S$.  Without loss of generality, $j_0 = 1$.  Since
$w: = u_{\beta v \zhat_\mu, 1}\, \beta v \zhat_\mu$ does not appear on
the left-hand side of
$(\dagger)$, it must cancel out on the right.  Observe that
$w$ does not equal any path of the form  $u_{\alpha p \zhat_r, j}\,
\alpha p \zhat_r$ with
$\alpha p \zhat_r \ne \beta v \zhat_\mu$, for $v \zhat_\mu$ is the
longest initial subpath of $w$ which belongs to $\S$.  Nor does $w$
coincide with one of the paths 
$u_{\beta v \zhat_\mu, j} \beta v \zhat_\mu$ for $j \ne 1$. 
Consequently, $w$ must be one of the $u_{\alpha p \zhat_r, j} q \zhat_s$
for a $\S$-critical path $\alpha p\zhat_r$, some $q \zhat_s \in
\S(\alpha p \zhat_r)$, and some $j$ with $b_{\alpha p \zhat_r, j} \ne
0$.  On the other hand, 
$$\len(q \zhat_s) \ge \len(\alpha p \zhat_r) \ge \len (\beta v \zhat_\mu)$$
by the definition of $\S(\alpha p \zhat_r)$ and our minimality
assumption, whence $q \zhat_s$ is longer than $v \zhat_\mu$.  But this is
absurd as, by construction,
$v \zhat_\mu$ is the longest initial subpath of $w$ which belongs to $\S$. 
This contradiction completes the linear independence argument. \qed
\enddemo

We next describe an ideal $\I(\S)$ of polynomials in $K[X]$ such that
$\psi(\grassS)$ is the vanishing locus of $\I(\S)$. 

\definition{3.14 Polynomials for the image $\psi(\grassS)$ in
$\AA^N$}  

Let $R$ be any finite generating set for the ideal $I \subseteq K\Gamma$
of relations for $\la$, viewed as a left ideal of $K\Gamma$.  Such a
generating set exists, since  all paths of length $L+1$ belong to $I$.  It
is, moreover, harmless to assume that
$R$ =
$\bigcup_{1 \le j \le n} R e_j$.  Consider the subset
$\widehat{R}$ of $\Phat \subseteq \P$, which is defined as follows:
$$\widehat{R} = \{ \rho \zhat_r \mid \rho \in R,\ 1 \le r \le t,\ \rho
e(r) = \rho\}.$$
Again, $e(r)$ is the primitive idempotent norming
$\zhat_r$, that is, $\zhat_r$ is a path of length zero in $\Phat$, which
starts and ends in $e(r)$.  Recall that, for any $i \in \{1, \dots, n\}$,
the number of indices $r$ with $e(r) = e_i$ equals the multiplicity $t_i$
of $S_i$ in $T$.  Hence, certain of the relations in
$R$ will ``fan out" to ``multiple incarnations" in 
$\widehat{R}$, while the elements $\rho \in R$ with
$\rho e_j = \rho$ and $S_j$ not a summand of $T$ will
not make an appearance in $\widehat{R}$.  

By Lemma 3.13, there exist, for each $\rho \zhat_r \in \widehat{R}$,
unique polynomials $\tau^{\rho \zhat_r}_{q \zhat_s} (X) \in K[X]$ with
the property that
$\rho \zhat_r \seq \sum_{q\, \zhat_s \in \S} \tau^{\rho \zhat_r}_{q
\zhat_s}(X)\, q \zhat_s$.
Define $\I(\S) \subset K[X]$ to be the ideal
generated by all the $\tau^{\rho\, \zhat_r}_{q
\zhat_s}$, where $\rho \zhat_r$ traces $\widehat{R}$. 

It is readily seen that the ideal $\I(\S)$ in $K[X]$ does not depend
on our specific choice of $R$, a fact which will also emerge as a
consequence of Proposition 3.15.  
\enddefinition

Essentially, our proof of Theorem 3.12 rests on the fact that
$\psi(\grassS)$ equals the vanishing set of the ideal $\I(\S)$ in
$\AA^N$.  As a by-product, we find that $\I(\S)$ in fact plays a universal
role relative to path algebras based on $\Gamma$ and $I$, with
coefficients in an arbitrary commutative $K$-algebra. 

Indeed, our setup can be generalized as follows:  Letting $\A$ be any
commutative $K$-algebra, we consider the path algebra (resp., path
algebra modulo relations)
$\A \Gamma \cong \A \otimes_K K\Gamma$ (resp., $\A \Gamma/ \A I
\cong \A \otimes_K \la$) with coefficients in $\A$.  The
projective $\A \Gamma$-module $\Phat(\A) = \bigoplus_{1 \le r \le
t} \A \Gamma \,\zhat_r$, where the $\zhat_r$ are again top elements
normed by the vertices $e(r)$ in the quiver $\Gamma$, is free as an
$\A$-module, having as basis the set of all paths in $\Phat$ (as
introduced ahead of Definition 3.1).  It specializes to
$\P$ for $\A = K[X]$, and  to $\Phat$ for $\A = K$.  In particular, $\Phat
= \Phat(K) \subseteq \Phat(\A)$, so that $\S$ can also be considered as a
subset of $\Phat(\A)$.  Finally, given an element $\fraka =
(\fraka_\nu)_{\nu \in N} \in \A^N$, we let $U(\fraka)$ be the submodule of
$\Phat(\A)$ generated by the differences
$$\alpha p \zhat_r - \sum_{q \zhat_s \in \S(\alpha p \zhat_r)}
\fraka_{\alpha p \zhat_r, q \zhat_s}\, q \zhat_s,$$
where $\alpha p \zhat_r$ runs through the $\S$-critical paths in $\Phat$,
identified with the corresponding elements in $\Phat(\A)$.

\proclaim{Proposition 3.15}  The $K$-algebra $\A_0 = K[X]/\I(\S)$ has the
following universal property relative to the skeleton $\S$.  Given any
$K$-algebra homomorphism
$\eta$ from $\A_0$ to another commutative $K$-algebra $\A$, set $\fraka =
\bigl( \eta(\overline{X}_\nu) \bigr) \in \A^N$.  Then the factor module 
$$M = \Phat(\A) / U(\fraka)$$
is an $\A \otimes_K \la$-module with the following property: For $0
\le l \le L$, the layer  $\Jhat ^l M / \Jhat ^{l+1} M$ 
is a free $\A$-module having as
basis the residue classes of the $p \zhat_r$ in $\S_l$; here $\S_l$ is
the set of paths of length $l$ in $\S$ as before, and $\Jhat$ denotes
the ideal $\A \otimes_K J$ in $\A \otimes_K \la$.
\smallskip

Conversely:  If $\A$ is a commutative $K$-algebra and $M = \Phat(\A)/U$ is
an $\A \otimes_K \la$-module satisfying the above layer condition, then
there exists a unique $K$-algebra homomorphism $\eta: \A_0 \rightarrow \A$
such that $U = U\bigl( \eta(\overline{X}_\nu) \bigr)$.
\endproclaim

\demo{Proof} Let $\eta$ and $\fraka$ be as in
the first assertion.  Then $\fraka = (\fraka_\nu)_{\nu \in N}$ belongs to
the vanishing set  $V_\A\bigl(\I(\S) \bigr)$ of $\I(\S)$ in $\A^N$, i.e.,
$\tau(\fraka) = 0$ for all $\tau \in \I(\S)$.  The definition of $\I(\S)$
entails that $M = \Phat(\A)/U(\fraka)$ is annihilated by $I$, which makes
$M$ an $\A \otimes_K \la$-module.  In particular, $\Jhat^{L+1} M = 0$.

In a first step, we show that $M$ is a free $\A$-module having basis $p
\zhat_r + U(\fraka)$, where $p \zhat_r$ runs through $\S$.  Preceding
$\eta$ by the quotient map $K[X] \rightarrow \A_0$ yields a
$K$-algebra homomorphism $\overline{\eta}: K[X] \rightarrow \A$, which in
turn induces an evaluation map $\theta: \P = \Phat(K[X]) \rightarrow
\Phat(\A)$ from the free $K[X]$-module $\P$ to the free $\A$-module
$\Phat(\A)$.  It is semilinear in the sense that $\theta(\tau y) =
\overline{\eta}(\tau) \theta(y)$ for $y \in \P$ and $\tau \in K[X]$. 
When $\A$ is endowed with the $K[X]$-module structure induced by
$\overline{\eta}$, the map 
$\theta$ gives rise to an isomorphism of free $\A$-modules
$$\A \otimes_{K[X]} \P \rightarrow \Phat(\A), \ \ 1\otimes q\zhat_s
\mapsto q \zhat_s,$$  
which maps the preferred path basis of $\A
\otimes_{K[X]} \P$ to that of $\Phat(\A)$; here $q \zhat_s$ traces the
paths in $\Phat$. By construction,
$\theta$ takes the defining generators for the submodule
$\C \subseteq \P$ to the defining generators for the submodule $U(\fraka)
\subseteq \Phat(\A)$, and by Lemma 3.13, $\P/ \C$ is still free over
$K[X]$, namely on basis $p \zhat_r + \C$, where $p \zhat_r$ traces
$\S$.  Thus we obtain an isomorphism $\A \otimes_{K[X]}
\bigl(\P/ \C \bigr) \cong \Phat(\A) / U(\fraka)$ sending $p \zhat_r + \C$
to $p \zhat_r + U(\fraka)$ for all $p \zhat_r \in \S$.  We infer that the
latter residue classes, $p \zhat_r + U(\fraka)$, indeed form a basis for
$\Phat(\A) / U(\fraka)$ over $\A$.     

To establish the first claim of the proposition, it
now suffices to prove that, for $0 \le l \le L$, the $p \zhat_r$
in $\sigma_l$ are $\A$-linearly independent modulo $\Jhat^{l+1}
M$.  Assuming this to be false, we let $l_0$ be minimal with respect to
failure.  Then $l_0 \ge 1$, and it is harmless to assume that $l_0 =
L$.  For otherwise, we may enlarge the ideal $I \subseteq K\Gamma$ so
that it contains the paths of length $l_0 + 1$ next to the original
relations, replace $\S$ by $\bigcup_{0 \le l \le l_0} \sigma_l$, adjust
the ideal $\I(\S) \subseteq K[X]$ according to 3.14, and replace $M$ by
$M/\Jhat^{l_0+1} M$.  In this modified setup,
$l_0$ is still minimal with respect to failure of our independence
condition.  So we only need to refute the assumption that the elements $p
\zhat_r + U(\fraka)$, $p \zhat_r \in \S_L$, are
$\A$-linearly independent in $M$.  But this was already established in the
first paragraph.

For the second claim, let $U$ be an $(\A \otimes_K K\Gamma)$-submodule of
$\Phat(\A)$ such that $M = \Phat(\A) / U$ is an $(\A \otimes_K
\la)$-module with the described layer property.  In particular,
$M$ is then a free $\A$-module with basis $p
\zhat_r + U$, where $p \zhat_r$ traces $\S$.  Letting $\alpha p \zhat_r$
be any $\S$-critical path of length $l$, our hypothesis on the layers of
$M$ yields coefficients $\frakb_{\alpha p \zhat_r, q \zhat_s} \in \A$, for
$q \zhat_s \in \S_l$, such that
$$\biggl(\alpha p \zhat_r \ - \ \sum_{q \zhat_s \in \S_{l}} \frakb_{\alpha
p \zhat_r, q \zhat_s}\, q \zhat_s \biggr) + U \in \Jhat^{l+1} M.$$
Clearly, we may assume that the above sum involves only paths
ending in the same vertex as $\alpha p \zhat_r$, meaning that all of the
$q \zhat_s$'s making an appearance belong to $\S(\alpha p \zhat_r)$.
An induction on $l \le L$ thus yields elements $\fraka_\nu \in \A$, for
$\nu \in N$, such that
$$\alpha p \zhat_r \ - \ \sum_{q \zhat_s \in \S(\alpha p \zhat_r)}
\fraka_{\alpha p \zhat_r, q \zhat_s}\, q \zhat_s\ \in \ U
\tag\dagger$$ 
for any $\S$-critical path $\alpha p \zhat_r$.
Uniqueness of the point $\fraka = (\fraka_\nu)_{\nu \in N} \in \A^N$
follows from freeness of $M$ over $\A$.  
This provides us with a unique $K$-algebra homomorphism $\etahat: K[X]
\rightarrow
\A$ such that $\bigl( \etahat(X_\nu) \bigr) = \fraka$.  By construction,
the map $\theta: \P = \Phat(K[X]) \rightarrow \Phat(\A)$ induced by
$\etahat$ sends the generators of the submodule $\C \subseteq \P$ to the
differences displayed in $\bigl(\dagger \bigr)$, so $\theta(\C)
\subseteq U$.  It now follows from Lemma 3.13 that $\theta(y) - \sum_{q
\zhat_s \in
\S} \tau_{q \zhat_s}^y (\fraka) q \zhat_s \in U$ for all $y \in \P$. 
Since $I \Phat(\A) \subseteq U$ (because $M$ is an $\A \otimes_K
\la$-module), this yields $\sum_{q \zhat_s \in
\S} \tau_{q \zhat_s}^{\rho \zhat_r} (\fraka) q \zhat_s \in U$ for all
$\rho \zhat_r \in \Rhat$.  Linear independence of the cosets $q \zhat_s
+ U$ thus implies $\fraka \in V_\A \bigl( \I(\S) \bigr)$.  Therefore
$\etahat$ induces a $K$-algebra homomorphism $\eta: \A_0 \rightarrow \A$
such that $\bigl( \eta(\overline{X_\nu}) \bigr) = \fraka$.  By
construction of $\fraka$, the submodule $U(\fraka)$ of
$\Phat(\A)$ is contained in $U$.  But as we saw in the proof
of the first claim, $\Phat(\A)/U(\fraka)$ is free of rank $|\S|$ over
$\A$; this rank coincides with the free rank of $\Phat(\A) /
U$ by hypothesis.  We conclude $U = U(\fraka)$ as required. 
Finally, uniquesness of $\eta$ holds by the uniqueness of $\fraka$.       
\qed \enddemo

\demo{Proof of Theorem 3.12}  Specializing Proposition 3.15 to $\A =
K$, we find that the image of the map $\psi$ equals the vanishing set
$V\bigl( \I(\S) \bigr)$ of the ideal $\I(\S)\subseteq K[X]$ in $\AA^N =
\AA^N(K)$, and that the inverse $\psi^{-1}$ sends any point $c \in V
\bigl( \I(\S) \bigr)$ to the submodule
$C = U(c)$ of $\Phat$.  That $\psi$ and $\psi^{-1}$ are morphisms of
varieties follows from an argument analogous to that of
\cite{\GeomII}, Theorem A, an early precursor of Theorem 3.12.  Thus
$\psi$ is indeed an isomorphism.

That the diagram in the statement of the theorem is commutative is now
clear.  Suppose that $K_0$ is a subfield of $K$ with the property
that elements generating $I \subseteq K\Gamma$ as a left ideal can be
found in $K_0 \Gamma$.  Since the construction of the
generators $\rho^{\alpha p \zhat_r}_{q \zhat_s}$ for $\I(\S)$ described in
3.14 may be carried out by applying Lemma 3.13 to the $K_0$-algebra $K_0
\Gamma/ (I \cap K_0 \Gamma)$, these polynomials may be assumed to belong
to $K_0[X]$.  As already mentioned, the claim concerning the algorithmic
nature is backed by the code in \cite{\codes}. \qed
\enddemo

\subhead 3.D. Affine coordinates for the subvarieties $\biggrassS$ of
$\grassbd$ \endsubhead

In this section, we outline how the definitions and methods of the
preceding section can be adapted to the case where the distinguished
projective $\la$-module considered is $\bp = \bigoplus_{1 \le r \le d}
\la z_r$ whose top $\bp / J \bp$ equals the semisimple left $\la$-module 
of dimension vector
$\bd$.  As before, $\bphat = \bigoplus_{1 \le r \le d}
K \Gamma \zhat_r$ is the corresponding projective $K\Gamma$-module, and
each top element $z_r$ of $\bp$ or $\zhat_r$ of $\bphat$ is
normed by the primitive idempotent $e(r)$.  Moreover, we assume that
$z_r$ coincides with the image of $\zhat_r$ under the canonical epimorphism
$\bphat \rightarrow \bp = \bphat/ I \bphat$.  

Again, we fix an an abstract skeleton $\S$ with dimension vector $\bd$
and top $T$, but this time, we choose $\S$ as a subset of $\bphat$.  Recall
that the dimension vector of $T$ is $\bt = (t_1, \dots ,t_n)$, whence $t :=
\sum_i t_i$ is the $K$-dimension of $T$.  Thus, it is clearly harmless to
assume that the paths of length zero in $\S$ are $\zhat_1, \dots,
\zhat_t$, and to identify the distinguished projective $K\Gamma$-module
$\Phat$ covering
$T$ with the direct summand $\bigoplus_{1 \le r \le t} K \Gamma
\zhat_r$ of $\bphat$.  Next we define $\S$-criticality; we do so by
carrying over the first of the equivalent descriptions of a
$\S$-critical path in Section 3.C to the present
situation.  We restate the definition for emphasis.

\definition{Definition 3.16 (pendant to Definition 3.9) and Comments}  A
{\it $\S$-critical path\/} is a path $q \zhat_r$ in $\bphat \setminus \S$,
of length at most $L$, with the property that every proper initial subpath
of $q \zhat_r$ belongs to $\S$. 

Note that we obtain two types of $\S$-critical paths in the
enlarged scenario:  The $\S$-critical paths in $\Phat$; these are
precisely the $\S$-critical paths of positive length having the form
$\alpha p \zhat_r$, where $\alpha$ is an arrow and $p \zhat_r \in \S$.
In addition, we have the $\S$-critical paths in the
complementary summand $\bigoplus_{t+1 \le r \le d} K\Gamma \zhat_r$; these
are precisely the paths
$\zhat_{t+1} \dots, \zhat_d$ of length zero in $\bphat
\setminus \S$.  Thus  every
$\S$-critical path is of the form $u \zhat_r$, where $u$ now is a
path of length $\ge 0$ in $K\Gamma$.

Given a $\S$-critical path $u \zhat_r$, its {\it $\S$-set\/},
$\S(u \zhat_r)$, consists of the paths $q \zhat_s \in \S$ of length $\ge
\len(u)$ which end in the same vertex as $u$.  So, in contrast to the
situation of Section 3.C, the set $\S(u \zhat_r)$ may contain paths of
length zero unless $d_i = t_i$ whenever $t_i \ne 0$.  In fact, each element
$\zhat_r \in \Phat$ with $t < r \le d$ is $\S$-critical, and $\S(
\zhat_r)$ contains all those candidates $\zhat_s$ among $\zhat_1, \dots,
\zhat_t$ for which $e(r) = e(s)$.     
\enddefinition

As in Section 3.C, we denote by $N$ the
set of all pairs 
$(u \zhat_r, q \zhat_s)$ such that the first entry is a
$\S$-critical path and the second entry belongs to the
corresponding $\S$-set $\S(u \zhat_r)$.  Since the critical paths of
positive length play a different role from those of length zero, we split
up the index set accordingly:  $N = N_1
\sqcup N_0$, where 
$$N_1 = \{\bigl(\alpha p \zhat_r, q \zhat_s \bigr) \mid \alpha p \zhat_r
\ \S\text{-critical of positive length},\ q \zhat_s \in \S(\alpha p
\zhat_r)\}$$ and
$$N_0 = \{\bigl(\zhat_r, q \zhat_s \bigr) \mid \zhat_r \ \S\text{-critical 
of length}\ 0,\ q \zhat_s \in \S(\zhat_r)\}.$$

Refer to Examples 3.5 for illustration. 
  
To obtain a family of scalars in $\AA^N$ pinning down the isomorphism type
of
$\bp/C$ for a point $C \in \biggrassS$, we observe as in 3.11:  For
any $\S$-critical path $u \zhat_r$, there exist unique scalars $c_{u
\zhat_r, q \zhat_s} \in K$ such that 
$$u z_r + C = \bigl(\sum_{q \zhat_s \in \S(u \zhat_r)} c_{u \zhat_r, q
\zhat_s} q z_s\bigr) + C$$ 
in $\bp/C$.  The family of scalars
$(c_\nu)_{\nu \in N}$ thus obtained, as $u \zhat_r$ traces the
$\S$-critical paths, determines
$\bp/C$ up to isomorphism; more strongly, it pins down $C$ (see Theorem
3.17 below).   As a consequence, we again obtain a morphism of varieties 
$$\bpsi: \biggrassS \rightarrow \AA^N, \qquad C \mapsto
\bigl(c_\nu \bigr)_{\nu
\in N}.$$

The first part of the following sibling of Theorem 3.12
relates the big varieties $\biggrassS$ to the small versions.  
   
\proclaim{Theorem 3.17 (pendant to Theorem 3.12)}  Let $\S\subseteq \bphat$
be a skeleton with top
$T$ and dimension vector $\bd$. 
\smallskip

\noindent {\rm (1)} The variety $\biggrassS$ is isomorphic to $\grassS
\times
\AA^{N_0}$.  In particular, $\biggrassS$ is an affine variety that
differs from $\grassS$ only by a direct factor which is a full affine
space.
\smallskip  

\noindent {\rm (2)} The map $\bpsi: \biggrassS \rightarrow
\AA^{N}$ induces an isomorphism from $\biggrassS$ onto a closed subvariety
of $\AA^N$. Its inverse $\bpsi^{-1}$ sends any point $c =
\bigl(c_\nu \bigr)_{\nu \in N}$ to the submodule $U(c) \subseteq
\bp$ which is generated by the differences
$$u z_r -  \sum_{q\zhat_s \in \S(u \zhat_r)} c_{u \zhat_r, q \zhat_s} q
z_s,$$   where  $u \zhat_r$ runs through the $\S$-critical paths.  This
isomorphism makes the following diagram commutative:
$$\xymatrixrowsep{2pc}\xymatrixcolsep{4pc}
\xymatrix{
\Im(\bpsi) \save+<-6.6ex,0.2ex> \drop{\AA^N\supseteq} \restore 
\ar[dd]_{\bpsi^{-1}} \ar[dr]^-{\boldsymbol{\chi}}\\
 &\save+<23ex,0ex> \drop{\txt{ \{\rm isom.~types of
 $\la$-modules with skeleton $\S$\} }} \restore\\
\biggrassS \save+<-12.6ex,0.1ex> \drop{\grassbd\supseteq} \restore 
\ar[ur]_-{\phi} }$$  
\noindent  Here $\boldsymbol{\chi}$ is the map which sends any point $c$ in
$\Im(\psi)$ to the isomorphism class of $\bp/U(c)$, and $\phi$ the map
which sends any point $C \in \biggrassS$ to the class of $\bp / C$.
\smallskip

\noindent {\rm (3)} If $K_0$ is a subfield of $K$ such that generators for
the left ideal
$I$ can be found in $K_0 \Gamma$, then polynomials defining
$\bpsi(\biggrassS)$ in $\AA^N$ can be chosen in $K_0[X]$.
\endproclaim

In light of Theorem 3.12, this ``big" version of the affine cover
theorem is immediate from the following lemma.  Since the proof of the
latter is straightforward, we leave it to the reader.  

\proclaim{Lemma 3.18}  We let $\S$ be a skeleton as in {\rm Theorem
3.17}.  Retaining the preceding notation, we suppose that
$\zhat_1, \dots, \zhat_t$ are the paths of length zero in $\S$ and view
the direct summand $P = \bigoplus_{1 \le r \le t} \la z_r$ of $\bp$ as the
distinguished $\la$-projective cover of $T$.  

\noindent {\rm (1)} Each point $C
\in \biggrassS$ has a $\la$-module decomposition $C = C_1 \oplus
C_0$, where
$$C_1 = \sum_{\alpha p \zhat_r \ \S\text{-critical of positive length}}
\la \biggl(\alpha p z_r - \sum_{q \zhat_s \in \S(\alpha p \zhat_r)}
c_{\alpha p \zhat_r, q \zhat_s} q z_s \biggr)$$
and  
$$C_0 = \sum_{t+1 \le r \le d} \la \biggl( z_r - \sum_{q
\zhat_s \in \S(\zhat_r)} c_{\zhat_r ,q \zhat_s} q z_s \biggr)$$
for some point $c\in \AA^N$.
In particular, $C = U(c)$.
\smallskip

\noindent {\rm (2)}  $C_1 = C \cap P$ belongs to $\grassS$,
the affine subvariety of $\grasstd$ which is determined by the skeleton
$\S$, when the latter is viewed as a subset of $\Phat$.  Moreover, 
$C_0$ is isomorphic to $\bigoplus_{t+1 \le r \le d} \la z_r$, and hence is
a complement of $P$ in $\bp$.  Both $C_1$ and $C_0$ are uniquely determined
by $C$. 
\smallskip

\noindent{\rm (3)}  Conversely, given any element $C_1 \in \grassS
\subseteq \grasstd$ and any element $(c_\nu)_{\nu \in N_0} \in \AA^{N_0}$,
let $C_0$ be the $\la$-submodule of $\bp$
which is generated by the differences 
$$z_r - \sum_{q \zhat_s \in \S(\zhat_r)} c_{\zhat_r, q \zhat_s} qz_s$$ 
with $(\zhat_r, q \zhat_s)$ running through $N_0$.  Then $\bp = P
\oplus C_0$, and $C_1 \oplus C_0$ belongs to $\biggrassS$.  \qed
\endproclaim

\head 4.  Degenerations in the the Grassmannian module varieties
\endhead

One of the main areas of application of the Grassmannian module varieties
lies in the fact that they permit a useful alternate perspective on
degenerations.  Recall that a {\it degeneration\/} of a module $M \in
\lamod$ of dimension vector $\bd$, which is represented by a point $x \in
\modlabd$ say, is any module represented by a point in the $GL_d$-orbit
closure of $x$ in $\modlabd$.  As is well known, the  relation $\bigl( M
\le M'
\iff M$ degenerates to $M'\bigr)$ defines a partial order on the set of
isomorphism types of modules of a fixed dimension vector.  A degeneration
$M'$ of
$M$ is called {\it top-stable\/} (resp., {\it layer-stable\/}) in case
$M'/JM' = M/JM$ (resp., $\SS(M') = \SS(M)$).  It is well-known (see
\cite{\Zwara}) that any point $x'$ in the closure of $GL_d.x$ in $\modlabd$
is connected to $x$ by way of a rational curve.  This statement can be
improved in the projective setting, due to the following result of
Koll\'ar (\cite{\degen}, Proposition 3.6):  Whenever $V$ is a 
unirational irreducible projective variety over $K$ (meaning that the
function field of $V$ embeds into a purely transcendental extension field
of $K$), any two points in $V$ are linked by a curve isomorphic to
$\PP^1$; in other words, given $C$ and $C'$ in $V$, there exists a
morphism $\phi:
\PP^1 \rightarrow V$, such that $\Im(\phi)$ contains both $C$ and
$C'$.  Observe that all our acting automorphism groups are rational
varieties, whence the orbits in the considered Grassmannians under the 
actions of the pertinent automorphism groups of projectives are all
unirational, as are their (projective) closures.  Combining Propositions
2.1 and 2.5 with curve connectedness of the orbit closures, we obtain:

\proclaim{Proposition 4.1}  Let $M \in \lamod$ have top $T$ and dimension
vector $\bd$. 
Moreover, let $\bp$ be the distinguished projective cover of the
semisimple module with dimension vector $\bd$, and $Q$ a projective module
with $T \le Q/JQ \le \bp / J\bp$, such that $Q$ is in turn equipped with a
distinguished sequence of top elements, as in {\rm Section 2.D}. 
\smallskip

\noindent {\rm (a)}  If $M = Q/C$ with $C \in [\grassbd]_Q$, the following
statements are equivalent:

$\bullet$  $M'$ is a  degeneration of $M$ with $M'/JM' \le Q/JQ$;

$\bullet$  $M' \cong Q/C'$, where $C'$ belongs to the
closure of $\Aut_\la(Q).C$ in $[\grassbd]_Q$;

$\bullet$  $M' \cong Q/C'$, where $C$ and $C'$ belong to the image
of a curve $\PP^1 \rightarrow \overline{\Aut_\la(Q).C}$, the latter
being the closure of $\Aut_\la(Q).C$ in $[\grassbd]_Q$.
\smallskip

\noindent {\rm (b)}  If $M = \bp/C$ with $C \in \grassbd$, the following
statements are equivalent:  

$\bullet$  $M'$ is a degeneration of $M$;
 
$\bullet$  $M' \cong \bp / C'$, where $C'$ belongs to the closure of
$\bigautlap.C$ in $\grassbd$;

$\bullet$  $M' \cong \bp/ C'$, where $C$ and $C'$ belong to the image
of a curve $\PP^1 \rightarrow \overline{\bigautlap.C}$, the latter
being the closure of $\bigautlap.C$ in $\grassbd$. \qed \endproclaim

In the special case where $Q$ coincides with the projective cover $P$ of
$T$, the equivalences under (a) have already been used to advantage to
explore top-stable degenerations of $M$, both theoretically and
computationally (see \cite{\degen}).

There is a special class of degenerations which is particularly
accessible by way of the Grassmannian varieties:

\definition{Definition 4.2}  A degeneration $M'\cong \bp/C'$ of
$M \cong \bp/C$ is called {\it unipotent\/} in case $C'$ belongs to the
closure of $\bigunirad.C$ in $\grassbd$. \enddefinition

In Proposition 4.4. and Corollary 4.5, we will see that, in exploring
unipotent degenerations, we may again cut down to smaller
settings, depending on the bound we place on the size of the tops. 

As a first application of Rosenlicht's result on unipotent group actions
(Theorem 2.14(2)), we find that no module has a proper layer-stable
degeneration which is unipotent.

\proclaim{Theorem 4.3} If $M, M' \in \lamod$ such that $M'$ is a proper
unipotent degeneration of $M$, then $\SS(M')$ properly dominates $\SS(M)$
in the sense of {\rm Definition 2.10}. In particular, $\SS(M') \ne \SS(M)$.
\endproclaim

\demo{Proof}  Suppose that $M \cong \bp/C$ and $M' \cong \bp/C'$ with $C,
C' \in \grassbd$.  If $\SS(M) = \SS$, then $\bigunirad.C$ is a closed
subvariety of $\biggrassSS$ by Theorem 2.14(2): Indeed, the $\biggrassS$,
for $\S$ compatible with $\SS$, form an affine open cover of
$\biggrassSS$ which is stable under the $\bigunirad$-action, whence the
intesection
$\biggrassS \cap \bigl( \bigunirad.C \bigr)$ is closed for each $\S$. 
Since
$C' \notin
\bigautlap.C$, and a fortiori $C' \notin \bigunirad.C$, we therefore
conclude $C' \notin \biggrassSS$.  Consequently, $C' \in
\biggrass{\SS'}$ for some semisimple sequence $\SS'$ which properly
dominates $\SS$ by Observation 2.11. \qed \enddemo

As a consequence of Theorem 4.3, we find that the layer-stable
degenerations of $M$ in (\cite{\degen}, Example 5.8) are not unipotent.  On
the other hand, usually, there is a plethora of unipotent degenerations
(cf\. Example 4.6).

We next apply Koll\'ar's result to the closure
of $\bigl(\Aut_\la(Q) \bigr)_u.C$ in $[\grassbd]_Q$. The geometric
structure of such an orbit  being well understood by Theorem 2.14(1),
curve connectedness takes on a particularly user-friendly form.  It will be
convenient to return to the limit notation for extensions of nonsingular
rational curves in an irreducible projective variety, which was used in
\cite{\degen}, Proposition 4.2 and Corollary 5.4:  Namely, if
$U$ is a dense open subset of $K = \AA^1 \subset \PP^1$,
and $\rho: U \rightarrow
\overline{\bigl(\Aut_\la(Q) \bigr)_u.C}$ a morphism, then $\rho$
extends uniquely to a curve  $\overline{\rho}$ defined on $U \cup
\{\infty\}$; hence, it is unambiguous to write $\lim_{\tau \rightarrow
\infty} \rho(\tau)$ for $\overline{\rho}(\infty)$.

\proclaim{Proposition 4.4}  Let $C \in [\grassbd]_Q$.  Moreover, for each
of the distinguished top elements $z_r$ of $Q$, pick finitely many paths
$p_{ri} z_{j(r,i)}$ in $JQ$ such that, modulo $\Hom_\la (Q,
C\cap JQ)$, each map in the space $\Hom_\la(Q,JQ)$ is of the
form $z_r \mapsto \sum_i a_{ri} p_{ri} z_{j(r,i)}$ for suitable $a_{ri}
\in K$,  $1 \le r \le m$. {\rm{(}}In particular, for each $r$, we
only need to consider paths $p_{ri}$ that end in $e(r)$.{\rm{)}}   

Then the unipotent
degenerations of $Q/C$ with top $\le Q/JQ$ are precisely the quotients
$Q/C'$ with
$$C' \ = \ \lim_{\tau \rightarrow \infty} g_\tau(C), \ \text{where} \
g_\tau(z_r) = z_r + \sum_i a_{ri}(\tau) p_{ri} z_{j(r,i)};$$ 
here the $a_{ri}(\tau)$ trace the rational functions in $K(\tau)$.
\endproclaim 

\demo{Proof}  First we consider the case $Q = \bp$. Say $\bp/C'$ is a
unipotent degeneration of $\bp/C$. This means that both $C$ and $C'$
belong to the projective variety $\overline{\bigunirad.C}$, whence, by
Koll\'ar's result, there exists a
curve
$\phi:
\PP^1
\rightarrow
\overline{\bigunirad.C}$ such that $C'= \lim_{\tau\rightarrow\infty}
\phi(\tau)$. In light of the comments preceding Proposition 2.14, $\phi$
restricts to a curve $\tau \mapsto g_\tau(C)$ of the indicated format,
with $g_\tau(z_r) = z_r + \sum_i a_{ri}(\tau) p_{ri} z_{j(r,i)}$ on the
preimage of
$\bigunirad.C$ under $\phi$.  Hence,
the claim holds in this case.

Now let $Q$ be a direct summand of $\bp$.  Without loss of generality,
we may assume that the distinguished top elements $z_1, \dots, z_m$ of
$Q$ coincide with the first $m$ candidates on the distinguished list $z_1,
\dots, z_d$ for $\bp$.  In particular, $\bp = Q \bigoplus Q_1$, where $Q_1
= \bigoplus_{m+1 \le r \le d} \la z_r$.  As explained in the remarks
preceding Definition 2.7, we may identify the projective variety
$[\grassbd]_Q$ with an isomorphic copy inside $\grassbd$; namely, with the
closed subvariety of
$\grassbd$ consisting of those points $C$ which contain $Q_1$. In other
words, we identify the points in $[\grassbd]_Q$ with those points in
$\grassbd$ which have the form $C = D \oplus Q_1$.  That $M'$ be a
unipotent degeneration of
$M = Q/C$ with $C \in [\grassbd]_Q$ such that $M' / J M' \le Q / JQ$,
means that there exists a point $C'$ in the intersection of $[\grassbd]_Q$
with the closure of
$\bigunirad.C$ in $\grassbd = [\grassbd]_{\bp}$ such that $M'  \cong \bp /
C'$.  Write $C' = D' \oplus Q_1$.  Due to projectivity of the intersection,
we can find a curve  $$\psi: \PP^1 \rightarrow [\grassbd]_Q \cap
\overline{\bigunirad.C}$$ 
such that $C' = \lim_{\tau \rightarrow \infty} \psi(\tau)$.  By the
special case $Q = \bp$, the restriction of this curve to the preimage of
$\bigunirad.C$ under $\psi$ has the form $U: \tau \mapsto h_\tau(C)$, for
some dense open subset $U \subseteq K$ and an automorphism $h_\tau$ given
by $h_\tau(z_r) = z_r +
\sum_i a_{ri}(\tau) p_{ri} z_{j(r,i)}$ for $1 \le r \le d$, $1 \le j(r,i)
\le d$, suitable  paths $p_{ri} z_{j(r,i)}$ in $J\bp$ ending in $e(r)$,
respectively, and suitable $a_{ri}(\tau) \in K(\tau)$. We split
each $h_\tau$ into its components: $g_\tau \in \bigl( \Aut_\la(Q)
\bigr)_u$ plus
$h^{(1)}_\tau \in \bigl( \Aut_\la(Q_1) \bigr)_u$ plus $h^{(2)}_\tau \in
\Hom_\la(Q, JQ_1)$ plus
$h^{(3)}_\tau \in \Hom_\la(Q_1, JQ)$.  By the choice of $\psi$, we have
$Q_1 \subseteq h_\tau(C)$.  We infer $Q_1 \subseteq h_\tau(Q_1)$, and
hence $Q_1 = h_\tau(Q_1)$ for reasons of dimension.  Consequently,
$h^{(3)}_\tau = 0$.  Moreover, we do not alter the image $h_\tau(C)$ by
taking $h^{(1)}_\tau$ to be the identity on $Q_1$ and
$h^{(2)}_\tau$ to be zero.  It is therefore harmless to assume that $h_\tau
= g_\tau \oplus \id_{Q_1}$, which implies
$\lim_{\tau\rightarrow\infty} g_\tau(D) = D'$ for the unipotent
automorphisms
$g_\tau$ in $\Aut_\la(Q)$.  This completes the argument.               
\qed\enddemo

\proclaim{Corollary 4.5}  Suppose that $M \cong P/C$ with $C \in
\grasstd$ and that $M'$ is a top-stable unipotent degeneration of $M$. 
Then $M' = P/C'$, where $C'$ belongs to the closure of the orbit
$\unirad.C$ in $\grasstd$. \qed \endproclaim 

We close with an example to illustrate how Proposition 4.4 can be
applied to finding unipotent degenerations.

\definition{Example 4.6} Let $\Gamma$ be the quiver with a
single vertex labeled $1$ and two loops, $\alpha$ and
$\beta$.  We consider the special biserial algebra $\la = K\Gamma/I$,
where $I \subseteq K\Gamma$ is the ideal generated by $\alpha^2$,
$\beta^2$, and all paths of length $3$.  We denote the regular left
$\la$-module by $\la e_1$ to emphasize the side.  Let $M = \la e_1$ and $Q
= \bigoplus_{1 \le r \le 3} \la z_r \cong (\la e_1)^3$.  We will discuss
the unipotent degenerations of $M$ with top $S_1^3$.  For this purpose, we
alternatively present $M$ in the form $Q/C$ with $C = \la z_1 \oplus \la
z_3$ in $[\biggrass_5(\la)]_Q$.  

For $a \in K^*$, we denote by $M_a$ the band module on the word
$\alpha \beta^{-1} \alpha
\beta^{-1}$ corresponding to the $2\times 2$ Jordan matrix with
eigenvalue $a$; namely, $M_a = (\la z_1 \oplus \la z_2)/ U$, where $U$ is
generated by $\beta z_1 - \alpha z_2$ and $\beta z_2 - a^2 \alpha z_1 +
2a\, \alpha z_2$. We first show that, for any choice of $a$, the direct sum
$M_a\oplus S_1$ arises as a unipotent degeneration of $M$. In fact, the
following is a sequence of successive unipotent degenerations:

$$\xymatrixrowsep{1.5pc}\xymatrixcolsep{0.3pc}
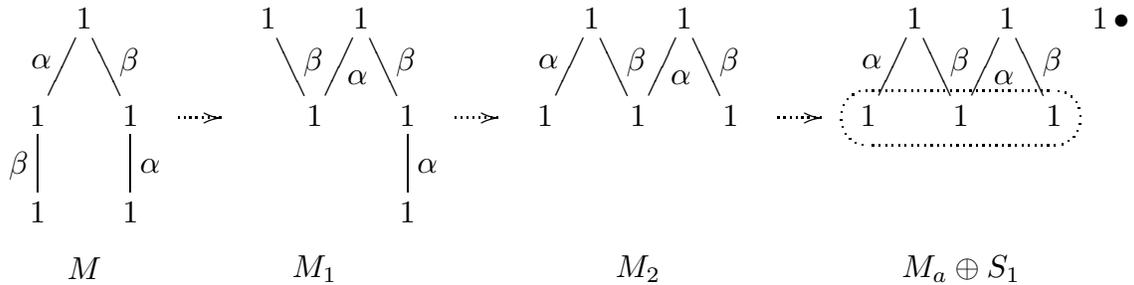
\xymatrix{
 &1 \edge[dl]_{\alpha} \edge[dr]^(0.55){\beta} &&&&&1 \edge[dr]^(0.55){\beta} &&1
\edge[dl]^{\alpha} \edge[dr]^(0.55){\beta} &&&&& &1 \edge[dl]_{\alpha}
\edge[dr]^(0.55){\beta} &&1 \edge[dl]^{\alpha} \edge[dr]^(0.55){\beta} 
&&&&& &1 \edge[dl]_{\alpha} \edge[dr]^(0.55){\beta} &&1 \edge[dl]^{\alpha}
\edge[dr]^(0.55){\beta} &&1\,\bullet \\
1 \edge[d]_{\beta} &&1 \edge[d]^{\alpha} &\dttdar[rr] &&&&1 &&1 \edge[d]^{\alpha}
&\dttdar[rr] &&&1 &&1 &&1 &\dttdar[rr] &&&1 \levelpool{4} &&1 &&1    \\
1 &\dropdown{M}&1 &&&&&\dropdown{M_1} &&1 &&&&&&\dropdown{M_2}
&&&&&&&& \dropdown{M_a\oplus S_1}
 }$$

Concerning the first link in the chain: $M_1 = Q/C_1$ with $C_1  =
\lim_{\tau \rightarrow \infty} g_\tau .C$, where, for $\tau \in K$, the map
$g_\tau
\in \bigl(\Aut_\la(Q) \bigr)_u$ is defined by $z_1
\mapsto z_1 + \tau(\beta z_1 - \alpha z_2)$ and $z_r \mapsto z_r$ for $r =
2,3$.   Indeed, from $z_1 \in C$, we first obtain $g_\tau(\alpha \beta z_1)
= \alpha \beta z_1
\in g_\tau .C$ and $g_\tau(\alpha z_1)= (\alpha + \tau \alpha \beta)z_1 \in
g_\tau .C$, whence $\alpha z_1 \in g_\tau .C$ for all $\tau$; this shows
$\alpha z_1
\in C_1$.  Then, using the computation technique for $1$-dimensional
subspaces of $\lim_{\tau \rightarrow \infty} g_\tau .C$ 
described in \cite{\degen}, Lemmas 4.7 and 5.5, we find 
$$\lim_{\tau
\rightarrow \infty} g_\tau(K z_1) = \lim_{\tau \rightarrow \infty}
K \bigl( \frac{1}{\tau} z_1 + (\beta z_1 - \alpha z_2) \bigr) = K(\beta z_1
- \alpha z_2) \subseteq C_1.$$
Thus, $C_1 \supseteq \la \alpha z_1 + \la (\beta z_1 -\alpha z_2)
+ \la z_3$. Since the latter sum is $10$-dimensional, and $Q/C_1$ is
$5$-dimensional, we infer equality. This shows $M_1$ to be a unipotent
degeneration of $M$.

To obtain $M_2$ as a unipotent degeneration of $M_1$, we
consider $g_\tau \in \bigl(\Aut_\la(Q) \bigr)_u$ for $\tau\in K$, defined
by $z_1 \mapsto z_1 +  \tau \beta z_2$, and $g_\tau(z_r) =
z_r$ for $r = 2,3$.  Moreover, we set $C_2 = \lim_{\tau \rightarrow
\infty} g_\tau.C_1$.  From $\alpha z_1 \in C_1$, we infer $\alpha z_1 +
\tau \alpha \beta z_2 \in g_{\tau}.C_1$, and hence, as above,
$$\lim_{\tau \rightarrow \infty} g_\tau(K \alpha z_1) = \lim_{\tau
\rightarrow \infty} K \bigl( \frac{1}{\tau} \alpha z_1 + \alpha \beta z_2
\bigr) = K \alpha \beta z_2 \subseteq C_2.$$
We next observe that $\beta z_1 - \alpha z_2 \in C_1$ implies
$\beta z_1 - \alpha z_2 = g_{\tau}(\beta z_1 - \alpha z_2)  \in
g_{\tau}.C_1$; moreover,
$\alpha z_1 \in C_1$ implies $\beta \alpha
z_1 = g_{\tau} (\beta \alpha z_1) \in g_{\tau}.C_1$ for all $\tau$. 
Consequently, these containments also hold in the limit, which shows
$$C_2 \supseteq \la \alpha \beta z_2 + \la(\beta z_1 - \alpha z_2) + \la
\beta \alpha z_1 + \la z_3.$$
A comparison of dimensions again gives equality.

To see that, for any $a \in K^*$, the module $M_2$ unipotently degenerates
to $M_a \oplus S_1$, we take $g_\tau(z_r) = z_r$ for $r = 1,2$ and
$g_\tau(z_3) = z_3 + \tau (\beta z_2 - a^2 \alpha z_1 +
2a\, \alpha z_2)$.  Due to $z_3\in C_2$, this yields 
$$\lim_{\tau \rightarrow \infty} g_\tau. C_2 = \la \beta \alpha z_1 + \la
(\beta z_1 - \alpha z_2) + \la (\beta z_2 - a^2 \alpha z_1 +
2a\, \alpha z_2) + \la \alpha z_3 + \la \beta z_3,$$
 as required.

Next we consider the $\PP^1$-family of modules $M_{[c_1{:}c_2]}$, which,
for
$c_1, c_2 \in K^*$, can be visulaized in the form
$$\xymatrixrowsep{1.5pc}\xymatrixcolsep{0.3pc}
\xymatrix{
1 \edge[dr]_{\alpha} &&1 \edge[dl]^{\beta} &&\dropdown{\oplus} &&&1
\edge@/_/[d]_{\alpha} \edge@/^/[d]^{\beta}  \\
 &1 &&&&&&1
}$$

\noindent While the modules $M_{[1{:}0]}$ and $M_{[0{:}1]}$, each a direct
sum of the displayed $3$-di\-men\-sion\-al string module and a string module of
dimension $2$, are unipotent degenerations of
$M$, we believe the remaining candidates not to be. 

All other $5$-dimensional left $\la$-modules with
$3$-dimensional top do arise as unipotent degenerations of $M$. There are
six sporadic ones, the $\PP^1$-family which we discussed at the outset
(next to the band modules $M_a \oplus S_1$ for $a \in K^*$, it contains two
string modules), and the following
$\PP^1\times\PP^1$-family:
$$\bigl( \la/\la (a_1\alpha-a_2\beta) \bigr) \oplus \bigl( \la/\la
(b_1\alpha-b_2\beta) \bigr) \oplus S_1\,, \qquad \bigl(\, [a_1{:}a_2],
[b_1{:}b_2] \in \PP^1 \,\bigr).  \tag\dagger$$
To deal with the latter, consider the family $(g_\tau)_{\tau\in K}$ in
$\bigl(\Aut_\la(Q) \bigr)_u$ with $g_\tau$ given by $z_1\mapsto z_1+ \tau
(a_1\alpha z_2- a_2\beta z_2)$, $z_2\mapsto z_2$, and $z_3 \mapsto z_3+
\tau (b_1\alpha z_1- b_2\beta z_1)$. Set $C'= \lim_{\tau
\rightarrow \infty} g_\tau.C$. As above, one obtains $a_1\alpha z_2-
a_2\beta z_2 \in C'$ and $b_1\alpha z_1- b_2\beta z_1 \in C'$ from $z_1,
z_3 \in C$. Moreover,
$\alpha z_3 \in C'$; indeed, $z_3\in C$ guarantees that $\alpha z_3- \tau
b_2 \alpha \beta z_1 \in g_\tau.C$ for all $\tau$. From
$z_1\in C$, we additionally derive $\alpha\beta z_1 \in g_\tau.C$ for
all $\tau$; this yields $\alpha z_3 \in g_\tau.C$ for all $\tau$, and hence
$\alpha z_3 \in C'$. Analogously, one shows that $\beta z_3\in C'$.
Consequently, 
$$C' \supseteq \la( a_1\alpha z_2-
a_2\beta z_2 ) +\la( b_1\alpha z_1- b_2\beta z_1 ) +\la \alpha z_3+
\la\beta z_3.$$
Since the right hand module is $10$-dimensional and $\dim Q/C'=5$, we
again derive equality, that is, $Q/C'$ equals a typical member of the
family $(\dagger)$.

Finally, we mention that, in particular, the modules $N_1$ and $N_2$ with
graphs
$$\xymatrixrowsep{1.5pc}\xymatrixcolsep{0.3pc}
\xymatrix{
1 \edge[dr]_{\alpha} &&1 \edge[dl]^(0.6){\beta} \edge[dr]^(0.4){\alpha} &&1
\edge[dl]^{\beta} &&&\dropdown{\txt{and}} &&&1 \edge[dr]_{\alpha} &&1
\edge[dl]^(0.6){\beta} \edge[dr]^{\alpha} &&1\,\bullet  \\
 &1 &&1 &&&&&&&&1 &&1
}$$

\noindent are unipotent degenerations of $M$. However, while $N_1$
degenerates to
$N_2$, this degeneration is not unipotent, by Theorem 4.3, since $\SS(N_1)=
\SS(N_2)$.            
\enddefinition

\enddocument